\documentclass[mn,fleqn,a4paper]{w-art}
\usepackage{times}
\usepackage{w-thm}

\usepackage[utf8]{inputenc}
\usepackage{latexsym}
\usepackage{amsthm}
\usepackage{amssymb}
\usepackage{amsmath}
\usepackage{amsfonts}
\newcommand{\tmsamp}[1]{\textit{#1}}
\theoremstyle{definition}
 \newtheorem{Def}[theorem]{Definition}
 \newtheorem{The}[theorem]{Theorem}
 \newtheorem{Pro}[theorem]{Proposition}
 \newtheorem{Lem}[theorem]{Lemma}
 \newtheorem{Cor}[theorem]{Corollary}
 \newtheorem{Rem}[theorem]{Remark}

\DeclareMathOperator{\diam}{diam}
\DeclareMathOperator{\supp}{supp}
\DeclareMathOperator{\Inv}{Inv}

\DeclareMathOperator{\cl}{cl}
\DeclareMathOperator{\A}{A}

\def\Rset{{\mathbb R}}

\def\Cset{{\mathbb C}}
\def\Nset{{\mathbb N}}
\def\Kset{{\mathbb K}}
\def\Qset{{\mathbb Q}}

\begin{document}

\title[Colombeau's algebra]{Generalized Solutions of a Nonlinear
  Parabolic Equation with\\ Generalized Functions as Initial Data}
\author[Sh. Jorge Aragona]{Jorge Aragona\footnote{Corresponding
     author: e-mail: {\sf aragona@ime.usp.br},}\inst{1}}
 \address[\inst{1}]{Universidade de São Paulo, Instituto de
   Matemática e Estatística}
\author[Sh. Antônio Ronaldo Gomes Garcia]{Antônio Ronaldo Gomes
  Garcia\footnote{Corresponding author: e-mail: {\sf ronaldogarcia@ufersa.edu.br}}\inst{2}}
\address[\inst{2}]{Universidade do Estado do Rio Grande do Norte,
  Departamento de Matemática e Estatística\\
Universidade Federal Rural do Semi-Árido, Departamento de Ciências Ambientais}
\author[Sh. Stanley Orlando Juriaans]{Stanley Orlando Juriaans
\footnote{Corresponding author: e-mail: {\sf ostanley@ime.usp.br},}\inst{1}}

\typeout{\dedicatory{This work that I am collaborate is dedicated my advisor 
Stanley Orlando Juriaans, your words, devotion and experience in
research in mathematics had been very important for my formation in
mathematics. I am grateful also the institution of promotion CNPQ-BRAZIL for 
the indispensable financial support during the elaboration of my 
doctorate thesis \cite{G}. We can not forget to thank those who gave us
the gift of life, the kind my God I am grateful.}}



\begin{abstract}
In \cite{bf} Br\'ezis and Friedman prove that certain nonlinear
parabolic equations, with the $\delta$-measure as initial data,  
have no solution. However in \cite{cl} Colombeau and Langlais prove 
that these equations have a unique solution even if the 
$\delta$-measure is substituted by any Colombeau generalized 
function of compact support. Here we generalize Colombeau and 
Langlais their result proving that we may take any generalized 
function as the initial data. Our approach relies on resent 
algebraic and topological developments of the theory of Colombeau 
generalized functions and results from \cite{A}. 
\end{abstract}

\maketitle

\section*{Introduction}  
\hspace{4.1mm} The necessity to prove the existence of solutions of 
equations may lead to the discovery of new and interesting
mathematical structures. Observe that, in general, it takes 
time for these structures to be fully understood and appreciated 
by the mathematical community. One might say that the algebras of 
generalized numbers and functions, introduced in the eighties by 
J.~F.~Colombeau, are among these structures. 
They are natural environments where the multiplication of 
distributions is something well defined and equals the classical 
definition for $\mathcal{C}^{\infty}$ functions. It should be noticed that 
regularization is definitely not the same as working in a Colombeau 
algebra. A Colombeau algebra is an algebraic and analytic environment! 

It is well known that rather simple and non-pathological linear 
equations have no distributional solution. However, the Colombeau algebras are 
environments where the concept
of derivation and that of solution of a P.D.E. can be generalized
in a natural way allowing, in many cases, to prove the existence of new
and interesting solutions for these equations (see \cite{arapan} and \cite{cl}).
  
This was proved by Colombeau and stimulated research in this 
new field as one can see from the results obtained by the 
Austrian, Brazilian, French, Serbian and South African research
groups and their collaborators (see \cite{ab} and \cite{kunz}).

In the beginning of the eighties Br\'ezis and Friedman showed that 
certain nonlinear parabolic equations have no solution if one chooses 
the initial data to be the $\delta$-measure. These non existence 
results, new in those days, were considered rather surprising because 
of several facts carefully explained by Br\'ezis and Friedman. An 
explanation for these non existence results was given in \cite{cl} 
by Colombeau and Langlais. Even more, they proved that the
Br\'ezis-Friedman equations do have a unique solution in the 
Colombeau algebra, as long as the initial data, which could be 
a distribution or Colombeau generalized function, had compact
support. One natural question can be formulated: Is there still a 
solution if the initial data has non-compact support? Is this solution 
unique too?

Since Colombeau defined his new algebras the Theory of Colombeau 
Generalized Functions has undergone rapid developments. The algebraic 
and topological aspects of the theory were developed which, in their 
turn allowed further development of other parts of the theory and 
led, to the development of a differential calculus which behaves like 
classical calculus (see \cite{difcal} and its reference list). 
A result obtained in \cite{ajos} states that the set of generalized 
functions of compact support is a dense ideal in the simplified
Colombeau algebra. 

To give an answer to the questions raised above we first generalize 
the  result mentioned in the last paragraph, i.e, we prove that the 
set of generalized functions of compact support is a dense ideal in 
the full Colombeau algebra.  Then, using the results of \cite{A} on 
quasi-regular sets, we push further the topological stepping stone 
of the theory. All this is done in  sections \ref{complet} and 
\ref{auxiliary}. In the last section (section \ref{ibvp}) we settle, 
in the positive, the two questions raised above. Notation is mostly 
standard unless explicitly stated.

\section{The completeness of the full Colombeau algebras}
\label{complet}
\hspace{4.1mm} As said in the introduction, in order to solve our
problem we will need to establish some topological facts about the 
Colombeau algebras on the closure of an open set. 

We begin by introducing a natural topology 
$\mathcal{T}_{\overline{\omega},b}$ on the algebra 
$\mathcal{G}\left(\overline{\omega}\right)$, where $\omega$ is a
\textit{bounded} open subset of $\Rset^{m}$. The basic
facts about the algebra $\mathcal{G}\left(\overline{\Omega}\right)$ are
presented in \cite{A} (for an arbitrary non-void open subset $\Omega$ of
$\Rset^{m}$).
The main results in this section are the completeness of
$\left(\mathcal{G}\left(\overline{\omega}\right),\mathcal{T}_{\overline{\omega},b}\right)$ and, as an easy consequence, the completeness of 
$\left(\mathcal{G}\left(\Omega\right),\mathcal{T}_{\Omega}\right)$, where 
$\mathcal{G}\left(\Omega\right)$ (resp. $\mathcal{T}_{\Omega}$)
is defined in [\cite{ab}, ~2.1] (resp. [\cite{A-F-J}, ~Definitions 
3.3 and 3.7 and Theorem 3.6]). In what follows, we will denote by
$\Omega$ (resp. $\omega$) a non void
open (resp. bounded open) subset of $\Rset^{m}$. We also set
$\textsc{I}:=]0,1]$,  
$\textsc{I}_{\eta}:=\left]0,\eta\right[, ~\eta\in\textsc{I}, ~\Kset$
denotes $\Rset$ or $\Cset, ~\left\Vert f\right\Vert:=\left\Vert f\right\Vert_{\omega}=\left\Vert f\right\Vert_{\overline{\omega}}=\sup\{\left\vert f\left(x\right)
\right\vert|x\in\omega\},~\forall ~ f\in
\mathcal{C}\left(\overline{\omega},\Kset\right)$. If $K\subset
X\subset\Rset^{m}$, the symbol $K\subset\subset X$ means that $K$ is a compact
subset of $X$.  

In the sequel we will use freely the partial order relation  $\leq$ on
$\overline{\Rset}$ introduced in [\cite{A-F-J}, Lemma 2.1 and Definition 2.2].
Let us recall that for a given
$\varphi\in\mathcal{D}\left(\Rset^{m};\Kset\right),\varphi\neq 0$, we define
$$i\left(\varphi\right):=\diam\supp\left(\varphi\right)$$
and it is easily seen that $i\left(\varphi_{\varepsilon}\right)=\varepsilon 
i\left(\varphi\right)$ for all $\varepsilon >0$. For every $r\in\Rset$ the
function [\cite{A-F-J}, ~Ex.2.3] 
$$\widehat{\alpha_{r}^{\bullet}}:\varphi\in\A_{0}\longmapsto
i\left(\varphi\right)^{r}\in\Rset_{+}^{\ast}$$ is moderated and therefore 
$\alpha_{r}^{\bullet}:=\cl\left(\widehat{\alpha_{r}^{\bullet}}\right)\in
\overline{\Rset}_{+}^{\ast}$.

Since $\overline{\omega}\subset\subset\overline{\omega}$ the definitions
of $\mathcal{E}_{M}\left[\overline{\omega}\right]$ and
$\mathcal{N}\left[\overline{\omega}\right]$ (see [\cite{A}, definition
of $\mathcal{E}_{M}\left[X;\Kset\right]$ and Definition 1.2(c)]) becomes simpler: 
$u\in \mathcal{E}_{M}\left[\overline{\omega}\right]$ means  
{\tmsamp{
$$(\mbox{M})\qquad \left\vert \begin{array}{ll}
\forall~ \sigma \in\Nset^{m}~ \exists~ N\in\Nset~ \mbox{such that}~\forall
~\varphi \in \A_{N}~\exists ~ c=c\left(\varphi\right)>0~\mbox{and}\\
\eta=\eta\left(\varphi\right)\in\mathtt{I}~\mbox{verifying}~
\left\Vert\partial^{\sigma}u\left(\varphi_{\varepsilon },\cdot \right)\right\Vert \leq c\varepsilon^{-N}, ~\forall ~\varepsilon\in\mathtt{I}_{\eta};
\end{array}\right.$$
}} and $u\in\mathcal{N}\left[\overline{\omega}\right]$ means 
{\tmsamp{
$$(\mbox{N}) \qquad \left\vert \begin{array}{ll}\forall~ \sigma\in 
\Nset^{m} ~\exists ~N\in\Nset~ \mbox{and}~\gamma\in\Gamma ~\mbox{such that}, ~
\forall ~q\geq N ~\mbox{and}, ~\forall~ \varphi\in \A_{q} 
~\exists ~c=c\left(\varphi\right)>0\\
\mbox{and}~ \eta=\eta\left(\varphi\right)\in 
\text{\textsc{I} verifying }\left\Vert \partial
  ^{\sigma}u\left(\varphi_{\varepsilon},
\cdot\right)\right\Vert\leq c\varepsilon^{\gamma\left(q\right)-N},~\forall ~\text{ }\varepsilon \in \text{\textsc{I}}_{\eta} 
\end{array}\right.$$}}

In the remainder of this paper we shall use a fixed exhaustive sequence
$\left(\Omega_{l}\right)_{l\in\Nset}$ of open subsets of  $\Omega$ 
(see [\cite{A-F-J}, begin of section 2]).
\vspace{0.7cm}

{\bf{A natural topology on $\mathcal{G}\left(\overline{\omega}\right)$}}
\hspace{4.1mm} For a given
$f\in\mathcal{G}\left(\overline{\omega}\right)$, if 
$\widehat{f}\in\mathcal{E}_{M}\left[\overline{\omega}\right]$ 
is any representative of $f$ and $\sigma\in\Nset^{m}$ is arbitrary, we have
$$\partial^{\sigma}\widehat{f}\left(\varphi,\cdot\right)\in\mathcal{C}\left(\overline{\omega};\Kset\right),~\forall ~\varphi\in \A_{0},$$ hence
$$\left\Vert\partial^{\sigma}\widehat{f}\left(\varphi,\cdot\right)\right\Vert<+\infty,
~\forall ~\varphi\in \A_{0}.$$

\newpage

\begin{Lem}\label{complet-1.1} {\tmsamp{For every $u\in
  \mathcal{E}_{M}\left[\overline{\omega}\right]$ and $\sigma \in
  \Nset^{m}$,  
\begin{enumerate}
\item[$(a)$] The function $$u_{\sigma}:\varphi\in\A_{0}\longmapsto\left
\Vert\partial^{\sigma}u\left(\varphi,\cdot\right)\right\Vert
\in\Rset_{+}$$ is moderate (i.e., $u_{\sigma}\in\mathcal{E}_{M}
\left(\Rset\right)$ and so $\cl\left(u_{\sigma }\right)\in\overline{\Rset}_{+}$);
\item[$(b)$] If $v\in\mathcal{E}_{M}\left[\overline{\omega}\right]$ 
and $(u-v)\in\mathcal{N}\left[\overline{\omega}\right]$ then $(u_{\sigma}-v_{\sigma})\in\mathcal{N}\left(\Rset\right)$ (and so $\cl\left(u_{\sigma}\right)=\cl\left(v_{\sigma}\right)$).
\end{enumerate}}}
\end{Lem}
\begin{proof} Follows immediately from definitions and [\cite{A-F-J}, 
Lemma 3.1].
\end{proof}

\begin{Def}\label{complet-1.2} {\tmsamp{Fix any $\sigma\in\Nset^{m}$. For each 
$f\in\mathcal{G}\left(\overline{\omega}\right)$ we set 
$$\left\Vert f\right\Vert_{\sigma}:=\cl\left(\widehat{f}_{\sigma}\right)=\cl\left[\varphi\in \A_{0}\longmapsto\left\Vert\partial^{\sigma}\widehat{f}\left(\varphi,\cdot\right)\right\Vert\in\Rset_{+}\right]\in\overline{\Rset}_{+},$$
where $\widehat{f}$ is an arbitrary representative of $f$. For each $r\in\Rset$ we define 
$$W_{\sigma,r}=W_{\sigma,r}\left[0\right]:=\left\{f\in\mathcal{G}\left(\overline{\omega}\right)|\left\Vert
    f\right\Vert_{\kappa}\leq\alpha_{r}^{\bullet},~\forall~\kappa\leq\sigma\right\}.$$}} 
\end{Def}

\begin{Rem}\label{complet-1.3}
\begin{enumerate}
\item[$(a)$]{\tmsamp{If $\sigma\geq\sigma^{\prime}$
(i.e. $\sigma_{i}\geq\sigma_{i}^{\prime},~\forall~ i=1,2,\cdots,m$) then
$W_{\sigma,r}\subset W_{\sigma^{\prime},r}$ for each $r\in\Rset$; 
\item[$(b)$] If $r>r^{\prime}$ then $W_{\sigma,r}\subset W_{\sigma,r^{\prime}}$. 
So, in the filter basis $\left(W_{\sigma,r}\right)~ (\sigma\in\Nset^{m}$ and
$r\in\Rset$) we can replace the condition $r\in\Rset$ by $r\in\Nset^{\ast}$.}}
\end{enumerate}
\end{Rem}

In the sequel we will use the following trivial results about inequalities:
{\tmsamp{
\begin{enumerate}
\item[$(a)$] If $x,y\in\Rset$ then $x\leq y$ in $\Rset$ if and only if
  $x\leq y$ in $\overline{\Rset}$; 
\item[$(b)$] If \ $x,x_{1},y$ and $y_{1}$  belong to
  $\overline{\Rset}$ with $0\leq x\leq x_{1}$ and $0\leq y\leq y_{1}$ 
then $0\leq xy\leq x_{1}y_{1}$.
\end{enumerate}}}

\begin{Lem}\label{complet-1.4} {\tmsamp{For every $f,g\in
\mathcal{G}\left(\overline{\omega}\right)$ we have
\begin{enumerate}
\item[$(a)$] $\left\Vert f+g\right\Vert_{\sigma}\leq\left\Vert f\right\Vert_{\sigma}+\left\Vert g\right\Vert_{\sigma}$;
\item[$(b)$] $\left\Vert fg\right\Vert_{\sigma }\leq
\sum\limits_{\kappa\leq\sigma}\binom{\sigma}{\kappa}\left\Vert
f\right\Vert_{\kappa}\left\Vert g\right\Vert_{\sigma-\kappa}$.
\end{enumerate}}}
\end{Lem}

\begin{proof} $(a)$ From Definition \ref{complet-1.2} we have 
$\left\Vert f+g\right\Vert_{\sigma}=\cl\left(v\right)$, where
$$v:\varphi\in \A_{0}\longmapsto \left\Vert\partial^{\sigma }\left(\widehat{f}+\widehat{g}\right)\left(\varphi,\cdot\right)\right\Vert\in\Rset_{+}$$
and $\widehat{f}$ and $\widehat{g}$ are any representatives of $f$ and $g$
respectively. From Definition \ref{complet-1.2} we have also 
$\left\Vert f\right\Vert_{\sigma}+\left\Vert g\right\Vert_{\sigma
}=\cl\left(u\right)$, where $$u:\varphi\in
\A_{0}\longmapsto\left\Vert\partial^{\sigma}\widehat{f}
\left(\varphi,\cdot\right)\right\Vert +\left\Vert\partial^{\sigma}
\widehat{g}\left(\varphi,\cdot\right)\right\Vert\in\Rset_{+}.$$
Since $v\left(\varphi\right)\leq u\left(\varphi\right)$ for every 
$\varphi\in\A_{0}$, from [\cite{A-F-J}, Lemma 2.1 and Definition 2.2] 
 statement $(a)$ follows.

$(b)$ From Definition \ref{complet-1.2} it follows that 
$\left\Vert fg\right\Vert_{\sigma}=\cl\left(v\right)$, where
$$v:\varphi\in\A_{0}\longmapsto\left\Vert\partial^{\sigma}\left(\widehat{f}
\widehat{g}\right)\left(\varphi,\cdot\right)\right\Vert\in\Rset_{+}$$
($\widehat{f}$ and $\widehat{g}$ are arbitrary representatives of $f$ and
$g$ respectively) and that $\sum\limits_{\kappa\leq\sigma }\binom{\sigma
}{\kappa }\left\Vert f\right\Vert_{\kappa}\left\Vert
  g\right\Vert_{\sigma-\kappa}=\cl\left(u\right)$, where 
$$u:\varphi\in \A_{0}\longmapsto\sum\limits_{\kappa\leq\sigma
}\binom{\sigma}{\kappa}\left\Vert\partial^{\kappa}\widehat{f}
\left(\varphi,\cdot\right)\right\Vert\left\Vert\partial^{\sigma-\kappa}
\widehat{g}\left(\varphi,\cdot\right)\right\Vert\in\Rset_{+}.$$
Now, from Leibnitz formula (see [\cite{A}, [1.2], p.373]) we have 
$v\left(\varphi\right)\leq u\left(\varphi\right)$ for all 
$\varphi\in\A_{0}$ and therefore
$\cl\left(v\right)\leq\cl\left(u\right)$, which is statement $(b)$.
\end{proof}

The two following lemmas will be useful and the proofs, which are 
trivial (apply [\cite{A-F-J}, Lemma 2.1 and Example 2.3]), are 
omitted.

\begin{Lem}\label{complet-1.5}
\begin{enumerate}
\item[$(a)$]{\tmsamp{For every $k,r\in\Rset_{+}^{\ast}$ we have
$k\left(\alpha_{s}^{\bullet}\right)^{2}\leq\alpha_{r}^{\bullet}$
(resp. $k\alpha_{s}^{\bullet}\leq\alpha_{r}^{\bullet}$) whenever
$s>\frac{r}{2}$ (resp. $s>r$);
\item[$(b)$] For every $k,r\in\Rset_{+}^{\ast}$ and $N\in\Nset$ there
  is $s\in\Rset_{+}^{\ast}$ such that 
$k.\alpha_{-N}^{\bullet}.\alpha_{s}^{\bullet}\leq\alpha_{r}^{\bullet
}$ (it is suffices  to choose $s\geq N+r+1$).}}
\end{enumerate}
\end{Lem}

\begin{Lem}\label{complet-1.6} {\tmsamp{For every 
$g\in\mathcal{G}\left(\overline{\omega}\right)$ and $\sigma\in\Nset^{m}$
there are $c>0$ and $N\in\Nset$ such that $\left\Vert
  g\right\Vert_{\kappa}\leq c.\alpha_{-N}^{\bullet }$ for each 
$\kappa\leq\sigma$.}}
\end{Lem}

With the notation introduced in Definition \ref{complet-1.2} we have 
the following result.
\vspace{.4cm}
              
\begin{The}\label{complet-1.7} {\tmsamp{Let $\omega$ be an open bounded
subset of $\Rset^{m}$. Then the set
$$\mathcal{B}_{\overline{\omega},b}:=\{W_{\sigma ,r}|\sigma
\in\Nset^{m}~\mbox{and}~r\in \Nset^{\ast}\}$$ is a filter basis on
$\mathcal{G}\left(\overline{\omega}\right)$
which satisfies the seven conditions of}} \mbox{[\cite{A-F-J}, Proposition 1.2
(2$^{\underline{o}}$)]} {\tmsamp{and so determine a topology
$\mathcal{T}_{\overline{\omega},b}$
on $\mathcal{G}\left(\overline{\omega}\right)$ which is compatible with its 
$\overline{\Kset}$-algebra structure (here we assume that
$\overline{\Kset}$ is endowed with its topology}} 
(see [\cite{A-F-J}, Definition 2.10]) 
{\tmsamp{and $\mathcal{B}_{\overline{\omega },b}$ is a
fundamental system of $\mathcal{T}_{\overline{\omega},b}$-neighborhoods of
$0$ in $\mathcal{G}\left(\overline{\omega}\right)$. Moreover the topology 
$\mathcal{T}_{\overline{\omega },b}$ is metrizable}}
\footnote{{\tmsamp{Note that, in view of
Remark \ref{complet-1.3} $(b)$ the set
$\mathcal{B}_{\overline{\omega },b}^{\prime}:=\{W_{\sigma,r}|\sigma\in 
\Nset^{m} ~\mbox{and} ~r\in\Rset_{+}^{\ast}\}$
is another fundamental system of $\mathcal{T}_{\overline{\omega},b}$-neighborhoods of $0$ in $\mathcal{G}\left(\overline{\omega}\right)$.}}}.
\end{The}

In the proof below of Theorem \ref{complet-1.7}, for the sake of
simplicity, we write $\mathcal{B}$ instead of 
$\mathcal{B}_{\overline{\omega},b}$. Also we will use freely the
notation (GA$_{\text{\textsc{I}}}^{\prime}$),\dots,(AV$_{\text{\textsc{II}}}^{\prime
}$) for the seven condition in [\cite{A-F-J}, Proposition 1.2 
(2$^{\underline{0}}$)]. Note also that the proof below works for 
$\mathcal{B}_{\overline{\omega},b}^{\prime}$.
\begin{proof}
In view of [\cite{A-F-J}, Corollary 1.3] it is enough to show
that $\mathcal{B}$ is a filter basis which satisfies the four
conditions
(GA$_{\text{\textsc{I}}}^{\prime}$),~(GA$_{\text{\textsc{II}}}^{\prime
}$),~(AV$_{\text{\textsc{I}}}^{\prime }$) and
(AV$_{\text{\textsc{II}}}^{^{\prime }}$) of [\cite{A-F-J}, Proposition
1.2] and that the topology $\mathcal{T}_{\overline{\omega },b}$ determinate
by $\mathcal{B}$ induces on $\overline{\Kset}$ \ its own topology 
$\mathcal{T}$ (see [\cite{A-F-J}, Definition 2.10]).
That $\mathcal{B}$ is a filter basis follows at once since clearly $\mathcal{B\neq\varnothing}$ and $\varnothing\notin\mathcal{B}$ and, if
$W_{\sigma,r},W_{\tau,s}$ are any two elements of $\mathcal{B}$ then, by
defining $\kappa:=\max\left(\sigma,\tau\right)$ (i.e., 
$\kappa_{i}:=\max(\sigma_{i},\tau _{i}), \forall~i=1,2,...,m$) and 
$t:=\max\left(r,s\right)+1$, from
Remark \ref{complet-1.3} we have 
$W_{\kappa,t}\subset W_{\sigma,r}\cap W_{\tau,s}$.

\textbf{Verification of} (GA$_{\text{\textsc{I}}}^{\prime}$): Given any 
$W_{\sigma,r}\in\mathcal{B}$ it is enough to show that
for any $s>\frac{r}{2}$ we have $W_{\sigma,s}+W_{\sigma,s}\subset
W_{\sigma,r}$, which follows immediately
from Lemma \ref{complet-1.4} $(a)$ and Lemma \ref{complet-1.5} $(a)$.

\textbf{Verification of \ (}GA$_{\text{\textsc{II}}}^{\prime }$):
 Obvious since $U=-U$ for all $U\in\mathcal{B}$.

\textbf{Verification of \ }(AV$_{\text{\textsc{II}}}^{\prime}$)$:$ It
suffices to prove that $W_{\sigma,r}^{2}\subset W_{\sigma,r}$ for each
$W_{\sigma,r}\in\mathcal{B}$. Fixed any $W_{\sigma,r}\in\mathcal{B}$, 
for $\kappa\in\Nset^{m}$, let:
$$M_{\kappa}:=\max\limits_{r\leq\kappa}\binom{\kappa}{\tau};
~M:=\max\limits_{r\leq\kappa}M_{\kappa}; ~p_{\kappa}:=\mbox{number of terms
of}~\sum\limits_{r\leq\kappa}\binom{\kappa}{\tau}; ~p:=\max\limits_{\kappa\leq \sigma}p_{\kappa } ~\mbox{and}~ k:=M.p.$$
We must show that
\begin{equation}\label{complet-1.7.1}
f,g\in W_{\sigma,r}\Longrightarrow fg\in W_{\sigma,r}.
\end{equation} 

Indeed, the assumption in (\ref{complet-1.7.1}) means that 
$\left\Vert f\right\Vert_{\lambda}\leq\alpha_{r}^{\bullet}$ and 
$\left\Vert g\right\Vert_{\lambda}\leq\alpha_{r}^{\bullet }$ 
for all $\lambda\leq\sigma$ hence, from Lemma \ref{complet-1.4} $(b)$ 
and Lemma \ref{complet-1.5} $(a)$ [and the obvious fact:
$a,b\in\Rset_{+}^{\ast},a\leq b, ~\lambda,~\mu\in\overline{\Rset}_{+}$
and $\lambda\leq\mu\Longrightarrow a\lambda\leq b\mu$] we can conclude that
for every $\kappa\leq\sigma$  
$$\left\Vert fg\right\Vert_{\kappa}\leq\sum_{\tau\leq\kappa}\binom{\kappa}{\tau}
\left\Vert f\right\Vert_{\tau}\left\Vert g\right\Vert_{\kappa-\tau }\leq
p_{\kappa}M_{\kappa}\left(\alpha_{r}^{\bullet}\right)^{2}
\leq pM\left(\alpha_{r}^{\bullet}\right)^{2}\leq
k\left(\alpha_{r}^{\bullet}\right)^{2}\leq\alpha_{r}^{\bullet}$$ 
since $r>\frac{r}{2}.$ Therefore (\ref{complet-1.7.1}) holds.

\textbf{Verification of} (AV$_{\text{\textsc{I}}}^{\prime }$): For
given $g\in \mathcal{G}\left(\overline{\omega}\right)$ and $W_{\sigma,r}\in 
\mathcal{B}$, we will show that there exists 
$W_{\tau,s}\in\mathcal{B}$ such that $gW_{\tau,s}\subset
W_{\sigma,r}$. Indeed, from Lemma \ref{complet-1.6} it follows that we
can find $c>0$ and $N\in\Nset$ such that
\begin{equation}\label{complet-1.7.2}
\left\Vert g\right\Vert_{\tau}\leq c\alpha_{-N}^{\bullet}, 
~\forall ~\tau\leq\sigma.
\end{equation}

\noindent Let $M,~p,~M_{\kappa}$ and $p_{\kappa}$ be as in the preceding proof 
of (AV$_{\text{\textsc{II}}}^{\prime}$) and define
$k:=pMc$ ($c$ from (\ref{complet-1.7.2})). Now, associated to 
$k,~r$ (which appears in $W_{\sigma,r}$)
and $N\in \Nset$ (which appears in (\ref{complet-1.7.2})), from
Lemma \ref{complet-1.5} $(b)$ it follows that
we can find $s(:=N+r+1$, for instance) verifying
\begin{equation}\label{complet-1.7.3}
k.\alpha_{-N}^{\bullet}.\alpha_{N+r+1}^{\bullet}\leq\alpha_{r}^{\bullet}.
\end{equation}
Next, we define $W_{\tau,s}$ by $\tau:=\sigma$ and $s:=N+r+1$, that is,
$W_{\tau,s}:=W_{\sigma,N+r+1}$ and it remains to prove that
$gW_{\tau,s}\subset W_{\sigma,r}$. In fact, fix
any $f\in W_{\tau,s}=W_{\sigma,N+r+1}$: 
\begin{equation}\label{complet-1.7.4}
\left\Vert f\right\Vert_{\tau}\leq\alpha_{N+r+1}^{\bullet}, ~\forall
~\tau \leq \sigma
\end{equation}
 then, from Lemma \ref{complet-1.4} $(b)$, (\ref{complet-1.7.2}), 
(\ref{complet-1.7.4}), (\ref{complet-1.7.3}) and the definition
of $k$, we get for every $\kappa\leq\sigma$:
$$\left\Vert gf\right\Vert_{\kappa}\leq\sum\limits_{\tau\leq\kappa}\binom{\kappa}{\tau}\left\Vert g\right\Vert_{\tau}\left\Vert f\right\Vert_{\kappa-\tau}\leq p.M\left(c\alpha_{-N}^{\bullet}\right)\left(\alpha_{N+r+1}^{\bullet }\right)=k\alpha_{-N}^{\bullet}.\alpha_{N+r+1}^{\bullet}\leq\alpha_{r}^{\bullet},$$ hence $gf\in W_{\sigma,r}$ and therefore (AV$_{\text{\textsc{I}}}^{\prime}$) holds.

Since $\mathcal{B}$ satisfies (GA$_{\text{\textsc{I}}}^{\prime
}$),~(GA$_{\text{\textsc{II}}}^{\prime}$),~(AV$_{\text{\textsc{I}}}^{\prime}$)
and (AV$_{\text{\textsc{II}}}^{\prime}$) from [\cite{A-F-J}, Corollary 1.3], 
we know that $\mathcal{B}$ determines a topology $\mathcal{T}_{\overline{\omega },b}$ on $\mathcal{G}\left(\overline{\omega}\right)$ which is
compatible with the ring structure of
$\mathcal{G}\left(\overline{\omega}\right) $. Moreover, it is clear 
(see [\cite{A-F-J}, Definition 2.7]) that $W_{\sigma,r}\cap
\overline{\Kset}=V_{r}\left[0\right]$ for every
$\sigma\in\Nset^{m}$ and $r\in\Rset_{+}^{\ast}$, which
implies that the topology induced by
$\mathcal{T}_{\overline{\omega},b}$ 
on $\overline{\Kset}$ coincides with the
topology $\mathcal{T}$ (see [\cite{A-F-J}, Definition 2.10]). Therefore, once more
from [\cite{A-F-J}, Corollary 1.3], we can conclude that 
$\mathcal{T}_{\overline{\omega},b}$ is compatible with the structure
of $\overline{\Kset}$-algebra of $\mathcal{G}\left(\overline{\omega 
}\right)$, where $\overline{\Kset}$ is endowed with its own topology 
$\mathcal{T}$. The topology $\mathcal{T}_{\overline{\omega},b}$ 
is Hausdorff since obviously $$\left\Vert f\right\Vert_{\sigma}=0,
~\forall~\sigma\in\Nset^{m}\Longrightarrow f\equiv 0$$
and hence $\underset{\sigma,r}{\bigcap}W_{\sigma,r}=\{0\}$. Indeed, if $f\in
\underset{\sigma,r}{\bigcap}W_{\sigma,r}$ we have
$$\left\Vert f\right\Vert_{\sigma}\leq\alpha_{r}^{\bullet},
~\forall~\sigma,r\iff \left\Vert f\right\Vert_{\sigma}\in
V_{r}\left[0\right], ~\forall ~\sigma,r\iff ~\forall ~\sigma,
~\mbox{we have}~\left\Vert f\right\Vert_{\sigma}\in\bigcap_{r>0} V\left[0\right] =\{0\},$$ since $\mathcal{T}$ is Hausdorff.
Finally, it is clear that $\mathcal{B}_{\overline{\omega},b}$ and
$\mathcal{B}_{\overline{\omega},b}^{\prime}$ generate the same filter 
(of all the $\mathcal{T}_{\overline{\omega},b}$-neighborhoods of $0$) 
and since $\mathcal{B}_{\overline{\omega},b}$ is countable, it follows 
that (see [\cite{B}, chapter 9, $2^{\underline{0}}$ ed.,~section 1, 
$n^{\underline{0}}~ 4$, Proposition 2])
$\mathcal{T}_{\overline{\omega},b}$ is metrizable. 
\end{proof}

\begin{Lem}\label{complet-1.8} {\tmsamp{Let $\Omega$ and $\omega$ be
two open subsets of $\Rset^{m}$ such that 
$\overline{\omega}\subset\subset\Omega$. Then the restriction map 
$$r=r_{\overline{\omega}}^{\Omega}:f\in\mathcal{G}\left(\Omega\right)\mapsto 
f|_{\overline{\omega}}\in\mathcal{G}\left(\overline{\omega}\right)$$
is continuous with respect to the topologies $\mathcal{T}_{\Omega}$ and
$\mathcal{T}_{\overline{\omega},b}$ on $\mathcal{G}\left(\Omega\right)$ and
$\mathcal{G}\left(\overline{\omega}\right)$ respectively.}}
\end{Lem}

\begin{proof} Fix an arbitrary $W_{\sigma,r}\in\mathcal{B}_{\overline{\omega
    },b}$. It suffices  to show that for any $l\in\Nset$ such that 
$\overline{\omega}\subset\subset\Omega_{l}$ we have
$r\left(W_{l,r}^{\sigma}\right)\subset W_{\sigma,r}$, which follows at once from
the definitions (see the definition of $W_{l,r}^{\sigma }$ in 
[\cite{A-F-J}, Definition 3.3]).
\end{proof}

In the next results we denote by $_{\overline{\Kset}}$\textbf{Top
Alg} the category whose object are the $\overline{\Kset}$-topological
algebras with the natural morphisms.
\begin{Pro}\label{complet-1.9} {\tmsamp{The topology $\mathcal{T}_{\Omega }$ on 
$\mathcal{G}\left(\Omega\right)$ is the initial topology (in the category
$_{\overline{\Kset}}$\textbf{Top Alg}) for the family of homomorphisms 
$\left(r_{l}\right)_{l\in \Nset}$, where 
$$r_{l}=r_{\overline{\Omega}_{l}}^{\Omega}:f\in\mathcal{G}\left(\Omega\right)
\mapsto f|_{\overline{\Omega}_{l}}\in\mathcal{G}\left(\overline{\Omega}_{l}\right)$$
and $\mathcal{G}\left(\overline{\Omega}_{l}\right)$ is
endowed with the topology $\mathcal{T}_{\overline{\Omega}_{l},b}$. 
In other words, $\mathcal{T}_{\Omega }$ is the coarsest topology on 
$\mathcal{G}\left(\Omega\right)$, compatible with the
structure of $\overline{\Kset}$-algebra of $\mathcal{G}\left(\Omega
\right)$, for which all the maps $r_{l}$ ($l\in\Nset$) are continuous.}}
\end{Pro}
\begin{proof} Let denote by $\mathcal{T}_{\Omega}^{\ast}$ the initial
topology on $\mathcal{G}\left(\Omega\right)$ for the family
$\left(r_{l}\right)_{l\in\Nset}$ (in the category
$_{\overline{\Kset}}$\textbf{Top Alg.}) Then it is
well known that $\mathcal{T}_{\Omega}^{\ast}$ has a
fundamental system of $0$-neighborhoods $\mathcal{B}_{\Omega}^{\ast}$
 consisting of sets of the following type:
\begin{equation}\label{complet-1.9.1}
V=\bigcap\limits_{1\leq i\leq p}r_{l_{i}}^{-1}\left(W_{\sigma
    _{i}}^{\left(l_{i}\right)},_{r_{i}}\right) 
\end{equation}
where $p\in\Nset^{\ast},\left(l_{i}\right)_{1\leq i\leq
  p},~\left(\sigma_{i}\right)_{1\leq i\leq p}$ and $\left(
  r_{i}\right)_{1\leq i\leq p}$ are finite sequences in 
$\Nset,~\Nset^{m}$ and $\Nset$ respectively and $W_{\sigma,r}^{\left(l\right)}\in\mathcal{B}_{\overline{\Omega }_{l},b}$ (here we need the \textbf{upper index}
$\left(l\right)$ since we are working with the subset $W_{\sigma,r}$ of $\mathcal{G}\left(\overline{\Omega}_{l}\right)$ which is denoted by
$W_{\sigma,r}^{\left(l\right)}$). From Lemma \ref{complet-1.8} it
follows that all the maps $r_{l}\left(l\in\Nset\right)$ are continuous
when $\mathcal{G}\left(\Omega\right)$ is endowed with the topology
$\mathcal{T}_{\Omega}$ and therefore
$\mathcal{T}_{\Omega}^{\ast}\preccurlyeq \mathcal{T}_{\Omega}$. In
order to prove that $\mathcal{T}_{\Omega}\preccurlyeq\mathcal{T}_{\Omega }^{\ast }$ it suffices to show that
$$\forall ~W_{l,r}^{\sigma}\in\mathcal{B}_{\Omega}~\exists
~V\in\mathcal{B}_{\Omega}^{\ast}~\mbox{such that}~V\subset
W_{l,r}^{\sigma}.$$ 
In fact, we can prove the following more precise statement
\begin{equation}\label{complet-1.9.2}
r_{l}^{-1}\left(W_{\sigma,r}^{\left(l\right)}\right)=W_{l,r}^{\sigma},
~\left(\forall~\sigma,l,r\right);
\end{equation}
note that the first member of (\ref{complet-1.9.2}) is of the type
(\ref{complet-1.9.1}). Now it is clear that (\ref{complet-1.9.2})
follows immediately from the definitions.
\end{proof}

\begin{Lem}\label{complet-1.10} {\tmsamp{Let $\omega$ and $\omega_{1}$ 
be two bounded open subsets of $\Rset^{m}$ such that
$\overline{\omega}_{1}\subset\omega$.
Then the restriction map}} (see [\cite{A}, Definition 2.6] 
$$r_{\overline{\omega}_{1}}^{\overline{\omega}}:f\in\mathcal{G}
\left(\overline{\omega}\right)\mapsto f|_{\overline{\omega}_{1}}
\in\mathcal{G}\left(\overline{\omega}_{1}\right)$$ 
{\tmsamp{is continuous when $\mathcal{G}\left(\overline{\omega}\right)$ and 
$\mathcal{G}\left(\overline{\omega}_{1}\right)$ are endowed with the
topologies $T_{\overline{\omega},b}$ and $T_{\overline{\omega}_{1}},_{b}$ 
respectively.}}
\end{Lem}

\begin{proof} For given $W_{\sigma ,r}^{1}\in \mathcal{B}_{\overline{\omega
  }_{1}},_{b}$ one has that $r_{\overline{\omega
  }_{1}}^{\overline{\omega}}\left(W_{\sigma,r}\right)\subset W_{\sigma,r}^{1}$,
where $$W_{\sigma,r}:=\left\{f\in\mathcal{G}\left(\overline{\omega 
}\right)|\left\Vert f\right\Vert_{\tau}\leq\alpha_{r}^{\bullet}, ~\forall
  ~\tau\leq\sigma\right\}$$ belongs to  $\mathcal{B}_{\overline{\omega },b}$.
\end{proof}

The following result is an adaptation to our case of the category of the topological
metrizable $\Kset$-algebras of [\cite{H}, Chapter 2, section 11, Proposition 3].

\begin{Pro}\label{complet-1.11} {\tmsamp{Let $A$ be a metrizable
topological ring, $\left(F_{p}\right) _{p\in\Nset}$ a sequence
of metrizable topological $A$-algebras and
assume that the condition below holds:
\begin{enumerate}
\item[$(a)$] If $p,q\in\Nset$ and $p\leq q$ then there exists a continuous
homomorphisms of $A$-algebras $f_{pq}:F_{q}\rightarrow F_{p}$.\\
Now, for a given $A$-algebra $E$ we assume that the following two 
conditions hold:
\item[$(b)$] For each $p\in\Nset$ there exists an homomorphism of
$A$-algebras $f_{p}:E\rightarrow F_{p}$ such that $p,q\in\Nset$
and $p\leq q$ implies $f_{p}=f_{pq}\circ f_{q}$;
\item[$(c)$] If $\left(x_{p}\right)_{p\in\Nset}\in\prod\limits_{p\in\Nset}F_{p}$ 
and for each $p,q\in\Nset$ with $p\leq q$ we have
  $f_{pq}\left(x_{q}\right)=x_{p}$, 
then there is $x\in E$ such that $f_{p}\left(x\right)=x_{p}$ for every
  $p\in\Nset$. Moreover, assume that $E$ is endowed with the initial topology 
$\mathcal{T}$ for the sequence $\left(f_{p}\right)_{p\in\Nset}$
[in the category of the metrizable topological $A$-algebras]
and that $\mathcal{T}$ is metrizable.
\end{enumerate}
The if every $F_{p}$ is complete, $E$ is complete.}}
\end{Pro}

\begin{proof}
Let $\left(y_{\nu}\right)_{\nu\in\Nset}$ be a Cauchy sequence in
$E$ then by defining
\begin{equation}\label{complet-1.11.1}
y_{\nu }^{p}:=f_{p}\left(y_{\nu}\right)\in F_{p}
\end{equation}
it is clear that $\left(y_{\nu}^{p}\right)_{\nu\in\Nset}$ is a
Cauchy sequence in $F_{p}$ for each $p\in\Nset$ and hence
\begin{equation}\label{complet-1.11.2}
\exists ~x_{p}:=\lim\limits_{\nu\rightarrow\infty}y_{\nu}^{p}\in F_{p}.
\end{equation}
Now, we will prove that $x_{p}=f_{pq}\left(x_{q}\right)$ whenever $p,q\in 
\Nset$ and $p\leq q$. Fix $p,q\in\Nset$ with $p\leq q$ 
arbitrary. Since
$y_{\nu}^{q}\underset{\nu\rightarrow\infty}{\longrightarrow} x_{q}~\mbox{in}~ F_{q}$ the continuity of $f_{pq}$ shows  that
\begin{equation}\label{complet-1.11.3}
f_{pq}\left(y_{\nu}^{q}\right)\underset{\nu\rightarrow\infty}{\longrightarrow}
f_{pq}\left(x_{q}\right).
\end{equation} 
On the other hand, from (\ref{complet-1.11.1}) and the condition $(b)$ we get
\begin{equation}\label{complet-1.11.4}
f_{pq}\left(y_{\nu }^{q}\right)=f_{pq}\left(f_{q}\left(y_{\nu}\right)\right)=f_{p}\left(y_{\nu }\right)=y_{\nu}^{p}\underset{\nu
\rightarrow\infty}{\longrightarrow}x_{p}
\end{equation}
and then, $f_{pq}\left(x_{q}\right)=x_{p}$ follows from (\ref{complet-1.11.3}) and (\ref{complet-1.11.4}). Therefore, the
hypothesis $(c)$ implies that there exists $x\in E$ such that
\begin{equation}\label{complet-1.11.5}
f_{p}\left(x\right)=x_{p},~\forall~ p\in\Nset.
\end{equation}
Next, fix an arbitrary $\mathcal{T}$-neighborhood $W$ of $0$ in $E$ that we
can choose of the form
$$W:=\underset{1\leq k\leq n}{\bigcap} f_{p_{k}}^{-1}\left(U_{k}\right),$$ 
where $U_{k}$ is a $0$-neighborhood in $F_{p_{k}}\left(1\leq k\leq
  n\right)$. From (\ref{complet-1.11.2}) it follows
that for each $k=1,2,\dots, n$ there is $l_{k}\in\Nset$ so that
\begin{equation}\label{complet-1.11.6}
\nu \geq l_{k}\Longrightarrow \left(y_{\nu}^{p_{k}}-x_{p_{k}}\right)\in U_{k}
\end{equation}
and therefore, by defining $\nu_{0}:=\max\limits_{1\leq k\leq n}l_{k}$ ,
we can conclude from (\ref{complet-1.11.1}), (\ref{complet-1.11.5}) and (\ref{complet-1.11.6}) that
$$\nu\geq\nu_{0}\Longrightarrow (y_{\nu}-x)\in W,$$
hence $\lim\limits_{\nu\rightarrow\infty}y_{\nu }=x$.
\end{proof}

\begin{The}\label{complet-1.12} {\tmsamp{If $\omega$ is a bounded
open subset of $\Rset^{m}$ then
$\mathcal{G}\left(\overline{\omega}\right)$ 
endowed with the topology $\mathcal{T}_{\overline{\omega },b}$ is complete.}}
\end{The}
\vspace{.7cm}

Before the proof of Theorem~\ref{complet-1.12} we use it together with 
Proposition \ref{complet-1.11} to prove that 
$\mathcal{G}\left(\Omega\right)$ endowed with the topology 
$\mathcal{T}_{\Omega}$ is complete (where $\Omega$ is an arbitrary 
open subset of $\Rset^{m}$). The role of $A,~F_{p},~E,~f_{pq}$ and 
$f_{p\text{ }}$ will be performed here by $\overline{\Kset},~\mathcal{G}
\left(\overline{\Omega}_{l}\right),~\mathcal{G}\left(\Omega\right),
~r_{lj}:=r_{\overline{\Omega}_{l}}^{\overline{\Omega}_{j}}~\left(l\leq j\right)$ and 
$r_{l}:=r_{\overline{\Omega }_{l}}^{\Omega }~ \left(l\in\Nset\right)$, 
respectively. Hypotheses $(a)$ and $(b)$ of Proposition \ref{complet-1.11} are obviously satisfied in view of Lemma \ref{complet-1.10}
and, by Proposition \ref{complet-1.9}, we know that
$\mathcal{T}_{\Omega }$ is the initial topology on 
$\mathcal{G}\left(\Omega\right)$
for the sequence $\left(r_{l}\right)_{l\in\Nset}$ and 
$\mathcal{T}_{\Omega}$ is metrizable by Theorem \ref{complet-1.7}. 
So it remains to show that the hypothesis $(c)$ of Proposition 
\ref{complet-1.11} holds. Indeed, fix any 
$\left(g_{l}\right)_{l\in\Nset}\in\prod\limits_{l\in\Nset}\mathcal{G}
\left(\overline{\Omega}_{l}\right)$ verifying the condition 
$$l,j\in\Nset ~\mbox{and}~ l\leq j\Rightarrow
g_{j}|_{\overline{\Omega}_{l}}
=g_{l},$$ we will show that there exists $f\in \mathcal{G}\left(\Omega \right)$
such that $f|_{\overline{\Omega}_{l}}=g_{l}$ for each $l\in\Nset$.
Indeed, let $f_{l}:=g_{l}|_{\Omega_{l}}\in\mathcal{E}_{M}\left[\Omega_{l}\right]$, $\forall l\in\Nset,$ then it is clear that
$$l,j\in\Nset~\mbox{and} ~l\leq j\Rightarrow 
f_{j}|_{\Omega_{l}}=g_{j}|_{\Omega_{l}}=g_{l}|_{\Omega_{l}}=f_{l}.$$
Since $\mathcal{G}$ is a sheaf, there exists $f\in\mathcal{G}\left(\Omega
\right)$ such that $f|_{\Omega_{l}}=f_{l}, ~\forall ~l\in\Nset$
and therefore $$f|_{\Omega_{l+1}}=f_{l+1}=g_{l+1}|_{\Omega_{l+1}},$$
hence
$$f|_{\overline{\Omega}_{l}}=\left(g_{l+1}|_{\Omega_{l+1}}\right)|_{\overline{\Omega}_{l}}=g_{l+1}|_{\overline{\Omega}_{l}}=g_{l},
\forall~ l\in\Nset,$$ which shows that condition $(c)$ of Proposition
\ref{complet-1.11} holds in our case. 

We can then conclude that we have the following
consequence of Proposition \ref{complet-1.11} and Theorem \ref{complet-1.12}:
\begin{Cor}\label{complet-1.13} {\tmsamp{If $\Omega$ is an open subset of
$\Rset^{m}$ then $\mathcal{G}\left( \Omega \right)$ endowed with the topology 
$\mathcal{T}_{\Omega}$ is complete.}}
\end{Cor}

The remainder of this section is essentially devoted to the rather long proof
of Theorem \ref{complet-1.12} and to some auxiliary results which will need in 
the last section. 

In order to prove Theorem \ref{complet-1.12} we start from some 
easy remarks which will be useful later.
\begin{Rem}\label{complet-1.14}
\begin{enumerate}
\item[$(a)$]{\tmsamp{Let $J$ be an ideal of a commutative ring 
$A$ with identity and assume that $\Qset\subset A$. Then for each $a,b$ in $A/J$ and 
for every representative $\widehat{u}$ of $(a-b)$ there are
representatives $\widehat{a}$ and $\widehat{b}$ of $a$ and $b$
respectively such that $\widehat{u}=(\widehat{a}-\widehat{b})$.
\item[$(b)$] Let $\left(x_{n}\right)_{n\geqslant 1}$ be a Cauchy sequence in an
abelian topological metrizable group $G$ and $\left(
  W_{i}\right)_{i\geqslant 2}$ a sequence of $0$-neighborhoods in $G$
such that $W_{i+1}\subset W_{i}$ for all $i\geqslant 2$. Then there is
a subsequence $\left(x_{n_{i}}\right)_{i\geqslant 1}$ of 
$\left(x_{n}\right)_{n\geqslant 1}$ such that 
$\left(x_{n_{i+1}}-x_{n_{i}}\right)\in W_{i+1}$ for $i\geqslant 1$.}}
\end{enumerate}
\end{Rem}

In the result below, for $\varkappa\in\Nset^{m}$ we will use the following
function $\psi_{\varkappa}$ defined by
$\psi_{\varkappa}\left(x\right):=\left(\varkappa!\right)^{-1}x^{\varkappa},
~\forall ~x\in\Rset^{m}$ and, of course, the restriction $\psi _{\varkappa }|_{\overline{\omega}}$ 
will be also denoted by $\psi _{\varkappa}$. Note that
$\partial^{\varkappa}\psi_{\varkappa}\left(x\right)\equiv 1,~\forall~
x\in\Rset^{m}$ and
$\widehat{\theta}\psi_{\varkappa}\in\mathcal{N}\left[\overline{\omega}\right]$
for every $\widehat{\theta}\in\mathcal{N}\left(\Rset\right)$ where
$\left(\widehat{\theta}\psi_{\varkappa}\right)\left(\varphi,x\right):=\widehat{\theta}\left(\varphi\right)\psi_{\varkappa}\left(x\right), ~\forall~\left(\varphi,x\right)\in\A_{0}\times\overline{\omega }$.

In what follows we will often use the notation
introduced in Definition \ref{complet-1.2} and Lemma \ref{complet-1.1}.

\begin{Lem}\label{complet-1.15}{\tmsamp{Let $\omega$ be an open
      bounded subset of
$\Rset^{m},~g\in\mathcal{G}\left(\overline{\omega}\right)$ and fix an
arbitrary representative $\widehat{g}$ of $g$. Then the statements below holds:
\begin{enumerate}
\item[$(I)$] Assume that $\varkappa \in \Nset^{m}$ and that $\widehat{R}$ is
  any representative of $\left\Vert g\right\Vert_{\varkappa}$.
Then there are a representative $\widehat{g}^{\left(\varkappa\right)}$ of $g$
  and $\widehat{\theta}_{\varkappa}\in\mathcal{N}\left(\Rset\right)$ verifying the following conditions:
\begin{enumerate}
\item[$(a)$]
  $\left(\widehat{g}^{\left(\varkappa\right)}\right)_{\varkappa}=\widehat{R}$ 
(i.e.
  $\left\Vert\partial^{\varkappa}\widehat{g}^{\left(\varkappa\right)}\left(\varphi,\cdot\right)\right\Vert=\widehat{R}\left(\varphi\right),
  ~\forall ~\varphi \in \A_{0}$);
\item[$(b)$] $\widehat{g}^{\left(\varkappa\right)}=\widehat{g}+\widehat{\theta}_{\varkappa}\psi_{\varkappa}$; 
\item[$(c)$] $\left( \widehat{g}^{\left(\varkappa\right)}\right)_{\varkappa}=\widehat{g}_{\varkappa}+\widehat{\theta}_{\varkappa}$;
\item[$(d)$] If $g^{\ast\left(\varkappa\right)}:=\widehat{g}+\left(\widehat{\theta}_{\varkappa}-\left\vert\widehat{\theta}_{\varkappa
}\right\vert\right)\psi_{\varkappa}$, then
$\widehat{R}=\left(g^{\ast\left(\varkappa\right)}\right)_{\varkappa}+\left\vert\widehat{\theta}_{\varkappa}\right\vert$ 
\end{enumerate}
\item[$(II)$] For a given $\left(\sigma,r\right)\in\Nset^{m}\times\Nset^{\ast}$ the conditions below are equivalent:
\begin{enumerate}
\item[$(i)$] $\left\Vert g\right\Vert_{\varkappa }\leq\alpha_{r}^{\bullet}$;
\item[$(ii)$] There exist $\widehat{\theta}_{\varkappa}\in\mathcal{N}\left(\Rset\right)$ and a representative $g^{\ast\left(\varkappa\right) }$ of $g$ such that
\begin{equation}\label{complet-1.15.1}
\left\Vert\partial^{\varkappa}g^{\ast\left(\varkappa\right)}\left(\varphi,\cdot\right)\right\Vert
+\left\vert\widehat{\theta}_{\varkappa}\left(\varphi\right)\right\vert\leq\widehat{\alpha}_{r}^{\bullet}\left(\varphi\right),
~\forall ~\varphi \in \A_{0}
\end{equation}
and, moreover, $g^{\ast\left(\varkappa\right)}$ is given [as function of $\widehat{\theta}_{\varkappa}$ and
of the representative $\widehat{g}$ fixed at the beginning ] by $(I)$ $(d)$.
\end{enumerate}
\item[$(III)$] For a given $\left(\sigma,r\right)\in\Nset^{m}\times\Nset^{\ast}$ the conditions below are equivalent:
\begin{enumerate}
\item[$(i)$] $g\in W_{\sigma,r}$
\item[$(ii)$] There exist a finite sequence $\left(g^{\ast\left(\varkappa\right)}\right)_{\varkappa\leq\sigma }$ of
representatives of $g$ and a finite sequence
$\left(\widehat{\theta}_{\varkappa}\right)_{\varkappa\leq\sigma}$ in 
$\mathcal{N}\left(\Rset\right)$ such that
\begin{equation}\label{complet-1.15.2}
\left\Vert\partial^{\varkappa}g^{\ast\left(\varkappa\right)}\left(\varphi,\cdot\right)\right\Vert+\left\vert\widehat{\theta}_{\varkappa}\left(\varphi\right)\right\vert\leq\widehat{\alpha}_{r}^{\bullet}\left(\varphi\right),
~\forall~ \varphi\in \A_{0},~\forall~\varkappa \leq \sigma.
\end{equation}
\end{enumerate}
\end{enumerate}}}
\end{Lem}

Note that (\ref{complet-1.15.1}) means that
$\left(g^{\ast\left(\varkappa\right)}\right)_{\varkappa}+\left\vert\widehat{\theta}_{\varkappa}\right\vert\leq
\widehat{\alpha }_{r}^{\bullet }$ hence, from $(I)$ $(a)$, $(c)$ and $(d)$ it
is clear that (\ref{complet-1.15.1}) 
is equivalent to
\begin{equation}\label{complet-1.15.1.1} 
\left(\widehat{g}^{\left(\varkappa\right)}\right)_{\varkappa}=\widehat{g}_{\varkappa}+\widehat{\theta}_{\varkappa}\leq\widehat{\alpha}_{r}^{\bullet}.
\end{equation}
The notation $g^{\ast\left(\varkappa\right)},\widehat{g}^{\left(
\varkappa\right)}$ and $\widehat{\theta }_{\varkappa}$  emphasize that
 these functions depend on $\widehat{g}$ and $\varkappa$. 
Clearly (\ref{complet-1.15.2}) is equivalent to
\begin{equation}\label{complet-1.15.2.2}
\left(\widehat{g}^{\left(\varkappa\right)}\right)_{\varkappa}\left(\varphi\right)=\widehat{g}_{\varkappa}\left(\varphi
\right)+\widehat{\theta}_{\varkappa}\left(\varphi\right)\leq\widehat{\alpha}_{r}^{\bullet}\left(\varphi\right),
~\forall ~\varphi\in \A_{0}, ~\forall ~\varkappa \leq \sigma.
\end{equation}
\begin{proof}[\textbf{Proof of Lemma 1.15.}] $(I)$: Since $\widehat{g}$ is a
representative of $g$, from the definition of 
$\left\Vert g\right\Vert_{\varkappa}$ (see Definition \ref{complet-1.2})
it follows that there is a unique
$\widehat{\theta}_{\varkappa}\in \mathcal{N}\left(\Rset\right)$
such that
\begin{equation}\label{complet-1.15.3}
\widehat{R}=\widehat{g}_{\varkappa}+\widehat{\theta}_{\varkappa}
\end{equation}
Next, define $\widehat{g}^{\left(\varkappa\right)}$ by the identity in
$(b)$ then for each $\left(\varphi,x\right)\in \A_{0}\times\overline{\omega 
}$ we have
$$\partial^{\varkappa}\widehat{g}^{\left(\varkappa\right)}\left(\varphi,x\right)
=\partial^{\varkappa}\widehat{g}\left(\varphi,x\right)+
\widehat{\theta}_{\varkappa}\left(\varphi\right)$$
hence $$\left(\widehat{g}^{\left(\varkappa\right)}\right)_{\varkappa}
\left(\varphi\right)=\widehat{g}_{\varkappa}\left(\varphi\right)+
\widehat{\theta}_{\varkappa}\left(\varphi\right),
~\forall ~\varphi \in \A_{0}$$
which proves statement $(c)$ and, from (\ref{complet-1.15.3}), 
also statement $(a)$. Statement $(d)$ is also trivial since 
the definition of $g^{\ast\left(\varkappa\right)}$ shows that
$$g^{\ast\left(\varkappa\right)}+\left\vert\widehat{\theta}_{\varkappa}
\right\vert\psi_{\varkappa}=\widehat{g}+
\widehat{\theta}_{\varkappa}\psi_{\varkappa}$$ hence
$$\partial^{\varkappa}g^{\ast\left(\varkappa\right)}\left(
\varphi,x\right)+\left\vert\widehat{\theta}_{\varkappa}\left(\varphi
\right)\right\vert=\partial^{\varkappa}\widehat{g}\left(\varphi,x\right)+\widehat{\theta}_{\varkappa}\left(\varphi\right),
\forall ~\left(\varphi,x\right)\in \A_{0}\times\overline{\omega}$$ and so
$$\left(g^{\ast\left(\varkappa\right)}\right)_{\varkappa}+\left\vert 
\widehat{\theta}_{\varkappa}\right\vert=\widehat{g}_{\varkappa}+\widehat{\theta
}_{\varkappa }$$ and therefore $(d)$ follows from (\ref{complet-1.15.3}).

$(II)$: $(i)\Rightarrow (ii)$: Condition (i) means
(see [\cite{A-F-J}, Lemma 2.1 (iii)]) that
there exists a representative $\widehat{u}$ of 
$\left\Vert g\right\Vert_{\varkappa }-\alpha _{r}^{\bullet }$ such that
\begin{equation}\label{complet-1.15.4}
\widehat{u}\left(\varphi\right)\leq 0 ~\forall ~\varphi \in \A_{0}.
\end{equation}
From Remark \ref{complet-1.14} $(a)$ there are representatives 
$\widehat{R}_{1}$ and $\widehat{R}_{2}$ of $\left\Vert g\right\Vert_{\varkappa }$ and $\alpha_{r}^{\bullet}$
respectively such that
$\widehat{u}=\widehat{R}_{1}-\widehat{R}_{2}$. Since 
$\widehat{g}$ is a representative of $g$ it follows 
that $\widehat{g}_{\varkappa }$ is a representative of $\left\Vert
g\right\Vert _{\varkappa }$, hence
$$\widehat{n}_{\varkappa}:=\widehat{R}_{1}-\widehat{g}_{\varkappa}\in\mathcal{N}\left(\Rset\right)$$
which shows that $\widehat{u}=\left(\widehat{g}_{\varkappa }+\widehat{n}_{\varkappa}\right)-\widehat{R}_{2}$ and therefore from (\ref{complet-1.15.4}) we get
\begin{equation}\label{complet-1.15.5}
\widehat{u}\left(\varphi\right)
=\widehat{g}_{\varkappa}\left(\varphi\right)+\widehat{n}_{\varkappa
}\left(\varphi\right)-\widehat{R}_{2}\left(\varphi\right)\leq 0,
~\forall ~\varphi \in \A_{0}.
\end{equation}
Since $\widehat{h}:=\widehat{\alpha}_{r}^{\bullet}-\widehat{R}_{2}\in 
\mathcal{N}\left(\Rset\right)$ it follows that 
$\widehat{\theta}_{\varkappa}:=\widehat{h}+\widehat{n}_{\varkappa }\in\mathcal{N}\left(\Rset\right)$ hence
\begin{equation}\label{complet-1.15.6}
\widehat{R}:=\widehat{g}_{\varkappa }+\widehat{\theta}_{\varkappa}
\end{equation}
is a representative of $\left\Vert g\right\Vert_{\varkappa }$ and
therefore, the proof of $(I)$ implies that $\widehat{g}^{\left(\varkappa
\right)}:=\widehat{g}+\widehat{\theta}_{\varkappa}\psi_{\varkappa}$ (see $(I)$
$(b)$) satisfies
\begin{equation}\label{complet-1.15.7}
\widehat{R}=\left(\widehat{g}^{\left(\varkappa\right)}\right)_{\varkappa
}=\widehat{g}_{\varkappa }+\widehat{\theta}_{\varkappa }.
\end{equation}
Consequently 
\begin{eqnarray}
\left(\widehat{g}^{\left(\varkappa \right)}\right)_{\varkappa }-
\widehat{\alpha }_{r}^{\bullet }&=&\left(\widehat{g}_{\varkappa}+\widehat{
\theta }_{\varkappa}\right)-\left(\widehat{R}_{2}+\widehat{h}\right)\nonumber\\
&=&\left(\widehat{g}_{\varkappa }+\widehat{h}+\widehat{n}_{\varkappa}\right)
-\left(\widehat{R}_{2}+\widehat{h}\right)\nonumber\\
&=&\widehat{g}_{\varkappa}+\widehat{n}_{\varkappa}-\widehat{R}_{2}\nonumber
\end{eqnarray}
and therefore, from (\ref{complet-1.15.5}), we get
\begin{equation}\label{complet-1.15.8}
\left( \widehat{g}^{\left(\varkappa \right)}\right)_{\varkappa}\left(\varphi\right)
\leq\widehat{\alpha}_{r}^{\bullet}\left(\varphi\right), ~\forall~
\varphi \in \A_{0}.
\end{equation}
Finally, from $(I)$ $(d)$ we have $\widehat{R}=\left(g^{\ast\left(\varkappa
\right)}\right)_{\varkappa}+\left\vert\widehat{\theta}_{\varkappa}\right\vert$
which together with the 
first identity of (\ref{complet-1.15.7}) implies that
(\ref{complet-1.15.8}) become
$$\left(g^{\ast\left(\varkappa\right)}\right)_{\varkappa}
\left(\varphi\right)+\left\vert\widehat{\theta}_{\varkappa}\left(\varphi\right)\right\vert\leq\widehat{\alpha}_{r}^{\bullet}\left(\varphi\right),
\forall ~\varphi \in \A_{0},$$
which proves $(ii)$.

$(ii)\Rightarrow (i)$: Since $\left( g^{\ast \left(\varkappa \right)}\right)_{\varkappa }+\left\vert\widehat{\theta}_{\varkappa}\right\vert$ is a representative of $\left\Vert g\right\Vert_{\varkappa}$, by [\cite{A-F-J},
Lemma 2.1 (iii)] the statement (\ref{complet-1.15.1}) implies $(i)$.

$(III)$ From the definition of $W_{\sigma,r}$ condition $(i)$ is equivalent
to 
$$(i') \qquad\left\Vert g\right\Vert_{\varkappa}\leq\alpha_{r}^{\bullet}, ~\forall ~\varkappa\leq\sigma.$$
From $(II)$, for each fixed $\varkappa \leq \sigma$, the inequality
$$\left\Vert g\right\Vert_{\varkappa}\leq\alpha_{r}^{\bullet}$$
is equivalent to the statement
{\tmsamp{
$$"\exists
~\widehat{\theta}_{\varkappa}\in\mathcal{N}\left(\Rset\right),
\mbox{and} ~\exists ~\mbox{a representative} ~
g^{\ast\left(\varkappa\right)} ~\mbox{of} ~g~ \mbox{such that (\ref{complet-1.15.1})} ~\mbox{holds}",$$}}
hence $(i')$ is equivalent to $(III)~(ii)$.   
\end{proof} 

In what follows we will use the following notation
$$\sigma_{i}:=\left(i,i,...,i\right)\in\Nset^{m}, ~\forall ~i\in\Nset.$$

\begin{Lem}\label{complet-1.16} {\tmsamp{Let $\left(f_{\nu}\right)_{\nu\geq 1}$ be a
  Cauchy sequence in $\left(\mathcal{G}\left(\overline{\omega}\right), \mathcal{T}_{\overline{\omega },b}\right)$. Then there exists a
  subsequence $\left(f_{\nu_{i}}\right)_{i\geq 1}$ of $\left(f_{\nu}\right)_{\nu\geq 1}$ such that if $\widehat{f}_{\nu _{i}}$ is an arbitrary
representative of $f_{\nu_{i}}$ for each $i\geq 1$, then for each $i\geq 1$ there are a finite sequence $\left(R_{i}^{\ast\left(\varkappa\right)}\right)_{\varkappa\leq\sigma_{i+1}}$ of representatives of $f_{\nu_{i+1}}-f_{\nu_{i}}$ and
finite sequence $\left(\widehat{\theta}_{i\varkappa}\right)_{\varkappa\leq\sigma_{i+1}}$ in $\mathcal{N}\left(\Rset\right)$ such that
\begin{equation}\label{complet-1.16.1}
\widehat{f}_{\nu_{i+1}}-\widehat{f}_{\nu _{i}}=R_{i}^{\ast
\left(\varkappa\right)}+\left(\left\vert\widehat{\theta}_{i\varkappa
}\right\vert-\widehat{\theta}_{i\varkappa}\right)\psi_{\varkappa}, 
~\forall ~i\geq 1,~\forall ~\varkappa\leq\sigma _{i+1};
\end{equation}
\begin{equation}\label{complet-1.16.2}
\left\Vert\partial^{\varkappa}R_{i}^{\ast\left(\varkappa\right)}\left(\varphi,
\cdot\right)\right\Vert+\left\vert 
\widehat{\theta}_{i\varkappa}\left(\varphi\right)\right\vert\leq
\widehat{\alpha}_{i+1}^{\bullet}\left(\varphi\right),~\forall~
i\geq 1,~\forall, ~\varkappa\leq\sigma_{i+1}, ~\forall ~\varphi\in\A_{0};
\end{equation}
\begin{equation}\label{complet-1.16.3}
\left\Vert\partial^{\varkappa}\left(\widehat{f}_{\nu _{i+1}}-\widehat{f}_{\nu _{i}}\right)\left(\varphi,\cdot\right)
\right\Vert\leq 2.\widehat{\alpha}_{i+1}^{\bullet}\left(\varphi\right),
~\forall ~i\geq 1, ~\forall ~\varkappa\leq\sigma_{i+1},~\forall 
~\varphi\in\A_{0}.
\end{equation}}}
\end{Lem}
\begin{proof}
Since $\left(W_{\sigma_{i},{i}}\right)_{i\geq 1}$ is
a decreasing sequence of $0$-neighborhoods in
$\left(\mathcal{G}\left(\overline{\omega}\right),\mathcal{T}_{\overline{
\omega}},b\right)$, from Remark \ref{complet-1.14} $(b)$ it follows
that we can find a subsequence $\left(f_{\nu _{i}}\right)_{i\geq 1}$ of 
$\left(f_{\nu}\right)_{\nu\geq 1}$ such that
\begin{equation}\label{complet-1.16.4}
R_{i}:=\left(f_{\nu _{i+1}}-f_{\nu_{i}}\right)\in
W_{\sigma_{i+1},{i+1}}, ~\forall~i\geq 1.
\end{equation}
Fix an arbitrary representative $\widehat{f}_{\nu_{i}}$ of $f_{\nu_{i}}$
for each $i\geq 1$. We shall apply
Lemma \ref{complet-1.15} $(III)$ and (\ref{complet-1.16.4}) to the 
representative
$$\widehat{R}_{i}:=\widehat{f}_{\nu_{i+1}}-\widehat{f}_{\nu_{i}}$$
of $R_{i}$. Then we can find a finite sequence
$\left(\widehat{\theta}_{i\varkappa }\right)_{\varkappa\leq\sigma
  _{i+1}}$ in $\mathcal{N}\left(\Rset\right)$ and a finite
sequence
$\left(R_{i}^{\ast\left(\varkappa\right)}\right)_{\varkappa\leq\sigma_{i+1}}$
of representatives of $R_{i}$ such that
$$\left\Vert\partial^{\varkappa}R_{i}^{\ast\left(\varkappa\right)}\left(\varphi,\cdot\right)\right\Vert+\left\vert\widehat{\theta}_{i\varkappa}\left(\varphi\right)\right\vert\leq\widehat{\alpha
}_{i+1}^{\bullet }\left(\varphi\right), ~\forall ~i\geq 1, ~\forall~\varphi\in
\A_{0}, \mbox{and}~\forall~\varkappa\leq\sigma_{i+1},$$ which proves 
(\ref{complet-1.16.2}).
Moreover, by the proof of Lemma \ref{complet-1.15} $(III)$ 
we know that $R_{i}^{\ast\left(\varkappa\right)}$ is given by Lemma 
\ref{complet-1.15} $(I)$ $(d)$, that is
$$R_{i}^{\ast \left( \varkappa\right)}=\widehat{R}_{i}+\left(\widehat{\theta}_{i\varkappa}-\left\vert\widehat{\theta}_{i\varkappa}\right\vert\right)\psi_{\varkappa}\left(\varkappa\leq\sigma_{i+1,\text{ }}i\geq 1\right)$$
hence, from the definition of $\widehat{R}_{i}$ we get
$$\widehat{f}_{\nu_{i+1}}-\widehat{f}_{\nu_{i}}=R_{i}^{\ast\left(\varkappa\right)}+\left(\left\vert\widehat{\theta}_{i\varkappa }\right\vert-\widehat{\theta}_{i\varkappa}\right)\psi_{\varkappa}$$ 
which proves (\ref{complet-1.16.1}). Clearly from
(\ref{complet-1.16.2}) it follows that
\begin{equation}\label{complet-1.16.5}
\widehat{\alpha}_{i+1}^{\bullet}\left(\varphi\right)-\left\vert\widehat{\theta}_{i\varkappa}\left(\varphi\right)\right\vert\geq
0, ~\forall ~i\geq 1, ~\forall~\varkappa\leq\sigma_{i+1},~
\forall~\varphi\in \A_{0}
\end{equation}
which implies at once that
\begin{equation}\label{complet-1.16.6}
\widehat{\alpha}_{i+1}^{\bullet}\left(\varphi\right)-\widehat{\theta}_{i\varkappa
}\left(\varphi\right)\leq
2.\widehat{\alpha}_{i+1}^{\bullet}\left(\varphi\right),~\forall ~i\geq
1,~\forall ~\varkappa\leq\sigma_{i+1},~\forall~\varphi\in\A_{0}
\end{equation}
[in fact, otherwise (\ref{complet-1.16.6}) would be false 
implying at once a contradiction with (\ref{complet-1.16.5}).]
Now, from (\ref{complet-1.16.6}), ~(\ref{complet-1.16.1}) and (\ref{complet-1.16.2}) we can conclude that for all
$\left(\varphi,x\right)\in \A_{0}\times \overline{\omega},~i\geq 1$
and $\varkappa\leq\sigma_{i+1}$ we have
$$\partial^{\varkappa}\left(\widehat{f}_{\nu_{i+1}}-\widehat{f}_{\nu
_{i}}\right)\left(\varphi,x\right)=\partial^{\varkappa}R_{i}^{\ast
\left(\varkappa\right)}\left(\varphi,x\right)+\left(\left\vert 
\widehat{\theta}_{i\varkappa}\right\vert-\widehat{\theta}_{i\varkappa}\right) 
\left(\varphi\right),$$ hence
$$\left\Vert\partial^{\varkappa}\left(\widehat{f}_{\nu _{i+1}}-\widehat{f}_{\nu _{i}}\right)\left(\varphi,\cdot\right)\right\Vert+\widehat{\theta}_{i\varkappa}\left( \varphi \right) =\left\Vert \partial
^{\varkappa}R_{i}^{\ast\left(\varkappa\right)}\left(\varphi,\cdot\right)\right\Vert+\left\vert\widehat{\theta}_{i\varkappa}\left(\varphi\right)\right\vert\leq \ \widehat{\alpha}_{i+1}^{\bullet}\left(\varphi \right)$$
which from (\ref{complet-1.16.6}) implies
$$\left\Vert\partial^{\varkappa}\left(\widehat{f}_{\nu_{i+1}}-\widehat{f}_{\nu
    _{i}}\right)\left(\varphi,\cdot\right)\right\Vert\leq\widehat{\alpha}_{i+1}^{\bullet}\left(\varphi\right)-\widehat{\theta}_{i\varkappa}\left(\varphi\right)\leq
2.\widehat{\alpha }_{i+1}^{\bullet}\left(\varphi \right), ~\forall
~i\geq 1, ~\forall~\varkappa \leq \sigma _{i+1} ~\mbox{and}~\forall~\varphi \in \A_{0}$$
\end{proof}

The application of Lemma \ref{complet-1.16} in the proof of Theorem 
\ref{complet-1.12} shall show that
in fact, the important statement of Lemma \ref{complet-1.16} 
is (\ref{complet-1.16.3}), and that the
other statements (\ref{complet-1.16.1}) and (\ref{complet-1.16.2}) 
are only preparatory results.

\begin{proof}[\textbf{Proof of Theorem 1.12.}] Let $\left(f_{\nu}\right)_{\nu\geq
1}$ be a Cauchy sequence in
$\left(\mathcal{G}\left(\overline{\omega}\right),\mathcal{T}_{\overline{\omega},b}\right)$. 
 We must shows that there exists $f\in
\mathcal{G}\left(\overline{\omega}\right)$ such that 
$f_{\nu }\underset{\nu \rightarrow \infty}{\longrightarrow} f$, in the topology 
$\mathcal{T}_{\overline{\omega },b}$. To this end, we shall consider 
the subsequence $\left(f_{\nu_{i}}\right)_{i\geq 1}$ of 
$\left(f_{\nu}\right)_{\nu \geq 1}$ defined by (\ref{complet-1.16.4}) 
(see the proof of Lemma \ref{complet-1.16}) and fix a sequence 
$\left(\widehat{f}_{\nu _{i}}\right)_{i\geq 1}$ where 
$\widehat{f}_{\nu_{i}}$ is an arbitrary representative 
of $f_{\nu _{i}}$ for each $i\geq 1.$ We shall use the sequence 
$\left(\widehat{f}_{\nu_{i}}\right)_{i\geq 1}$
 to define an element
$\widehat{f}\in\mathcal{E}_{M}\left[\overline{\omega}\right]$ such
that $f_{\nu _{i}}\underset{i\rightarrow \infty}{\longrightarrow}
f:=\cl\left(\widehat{f}\right)\in\mathcal{G}\left(\overline{\omega}\right),$ 
 and hence $f_{\nu}\underset{\nu\rightarrow\infty}{\longrightarrow} f$.
So, from now on, we assume that the following are fixed: the subsequence $\left(f_{\nu_{i}}\right)
_{i\geq 1}$, a representative 
$\widehat{f}_{\nu_{i}}$ of $f_{\nu _{i}}$ for each $i\geq 1,
~R_{i}:=f_{\nu _{i+1}}-f_{\nu_{i}}$ and
$\widehat{R}_{i}:=\widehat{f}_{\nu_{i+1}}-\widehat{f}_{\nu_{i}}$. 
Since $\left(\nu_{i}\right)_{i\geq 1}$ is strictly increasing it is 
clear that $\left(\A_{\nu_{i}}\right)_{i\geq 1}$ and $\left(\text{\textsc{I}}_{\nu
_{i}^{-1}}\right)_{i\geq 1}$ are strictly decreasing. For every 
$i\geq 1$ and $x\in\overline{\omega}$ we define

$$\widehat{u}_{i}\left(\varphi,x\right):=\left\{ 
\begin{array}{cc}
\left(\widehat{f}_{\nu_{i+1}}-\widehat{f}_{\nu
    _{i}}\right)\left(\varphi,x\right),&
\mbox{if}~\left(\varphi,i\left(\varphi\right)\right)\in \A_{\nu_{i}}\times\text{\textsc{I}}_{\nu_{i}^{-1}}\\ 
 0, &\mbox{if} ~\left(\varphi,i\left(\varphi\right)\right)\notin \A_{\nu_{i}}\times \text{\textsc{I}}_{\nu _{i}^{-1}}\end{array}\right.$$
hence, for all $i\geq 1$, $x\in\overline{\omega }$ and $\beta
\in\Nset^{m}$ we have 
\begin{equation}\label{complet-1.12.1}
\partial^{\beta}\widehat{u}_{i}\left(\varphi,x\right)=\left\{ 
\begin{array}{cc}
\partial^{\beta}\widehat{R}_{i}\left(\varphi,x\right), & \mbox{if}~\left(\varphi ,i\left(\varphi\right)\right)\in \A_{\nu_{i}}\times\text{\textsc{I}}_{\nu _{i}^{-1}} \\ 
0, &\mbox{if} ~\left(\varphi,i\left(\varphi\right)\right)\notin
\A_{\nu_{i}}\times \text{\textsc{I}}_{\nu _{i}^{-1}}.
\end{array}\right.
\end{equation} 
Obviously, we have
$\widehat{u}_{i}\in\mathcal{E}_{M}\left[\overline{\omega}\right]$ for
all $i\geq 1$. Define
\begin{equation}\label{complet-1.12.2}
\widehat{f}\left(\varphi,x\right):=\widehat{f}_{\nu_{1}}\left(\varphi,x\right) 
+\sum\limits_{i\geq 1}\widehat{u}_{i}\left(\varphi,x\right), ~\forall
~\left(\varphi,x\right)\in \A_{0}\times\overline{\omega}.
\end{equation}
Clearly, for every $\left(\varphi,x\right)\in
\A_{0}\times\overline{\omega}$, 
the series in the second member of (\ref{complet-1.12.2})
is finite [indeed, note that for every $\varphi\in \A_{0}$ we have 
"or $\varphi\notin \A_{\nu_{1}} ~\mbox{or} ~\exists!~s\in\Nset^{\ast}$ 
such that $\varphi\in \A_{\nu_{s}}\cap\complement \A_{\nu_{s+1}}$"], 
hence $\widehat{f}$ is well defined and furthermore, 

\begin{equation} 
\label{complet-1.12.3}
\forall\beta
\in\Nset^{m}~\mbox{and}~\left(\varphi,x\right)\in
\A_{0}\times\overline{\omega} \mbox{\ we have that\ } 
\partial^{\beta}\left(\sum\limits_{i\geq 1}\widehat{u}_{i}\left(\varphi
,x\right)\right)=\sum\limits_{i\geq
1}\partial^{\beta}\widehat{u}_{i}\left(\varphi,x\right). 
\end{equation} 
Next we shall prove that
\begin{equation}\label{complet-1.12.4}
\widehat{f}\in \mathcal{E}_{M}\left[\overline{\omega}\right].
\end{equation}
Initially, note that from (\ref{complet-1.12.1}) we get
\begin{equation}\label{complet-1.12.5}
\left\Vert\partial^{\beta}\widehat{u}_{i}\left(\varphi,x\right)\right\Vert\leq
\left\Vert\partial^{\beta}\widehat{R}_{i}\left(\varphi,x\right)\right\Vert,
~\forall ~\varphi\in\A_{0},~\forall~ i\geq 1, ~\forall~\beta\in\Nset^{m}.
\end{equation}
Clearly, to prove (\ref{complet-1.12.4}) it suffices  to show the 
moderateness of the function
$$\widehat{U}=\widehat{U}_{1}:\left(\varphi,x\right)\in\A_{0}\times\overline{\omega}\longmapsto\sum\limits_{i\geq 1}\widehat{u}_{i}\left(\varphi,x\right)\in\Kset.$$
So, we must show the following statement {\tmsamp{
\begin{equation}\label{complet-1.12.6}
\left\vert\begin{array}{ll}
\mbox{For a fixed}~\beta\in\Nset^{m}, ~\exists~N\in\Nset~
\mbox{such that}~\forall ~\varphi\in\A_{N}, ~\exists ~ 
\widetilde{c}=\widetilde{c}\left(\varphi\right)>0 ~\mbox{and}\\
\exists ~\widetilde{\eta}=\widetilde{\eta}\left(\varphi\right)
\in \text{\textsc{I}} ~\mbox{verifying} 
~\left\Vert\partial^{\beta}\widehat{U}\left(\varphi_{\varepsilon },\cdot
\right)\right\Vert\leq\widetilde{c}\mathcal{\varepsilon}^{-N}, ~\forall
~\varepsilon\in \text{\textsc{I}}_{\widetilde{\eta
  }}.\end{array}\right.
\end{equation}}}
Indeed, for a given $\beta\in\Nset^{m}$, choose
$r\in\Nset^{\ast}$ 
such that
\begin{equation}\label{complet-1.12.7}
\beta\leq\sigma_{r+1}
\end{equation}
and consider the function $\widehat{U}_{r}$ defined by  
$\widehat{U}_{r}\left(\varphi,x\right):=\sum\limits_{i\geq
  r}\widehat{u}_{i}\left(\varphi,x\right), ~\forall
~\left(\varphi,x\right) 
\in\A_{0}\times\overline{\omega}$.
By the same finiteness argument used proving (\ref{complet-1.12.3}) we have
\begin{equation}\label{complet-1.12.8}
\partial^{\beta
}\widehat{U}_{r}\left(\varphi,x\right)=\sum\limits_{i\geq r}
\partial^{\beta}\widehat{u}_{i}\left(\varphi,x\right), ~\forall
~\left(\varphi,x\right)\in\A_{0}\times\overline{\omega}
\end{equation}
In view of (\ref{complet-1.16.3}) (see Lemma \ref{complet-1.16})
applied to $\widehat{R}_{i}$ with $\varkappa=\beta\leq\alpha_{i+1},
~\forall ~i\geq r$, see (\ref{complet-1.12.7})) we can write
\begin{equation}\label{complet-1.12.9}
\left\Vert\partial^{\beta}\widehat{R}_{i}\left(\varphi,\cdot\right)\right\Vert
\leq 2\widehat{\alpha }_{i+1}^{\bullet}\left(\varphi\right), ~\forall
~i\geq r, ~\forall ~\varphi\in\A_{0}.
\end{equation}
From (\ref{complet-1.12.8}) we have
$$\left\vert\partial^{\beta}\widehat{U}_{r}\left(\varphi,x\right)\right\vert\leq\sum\limits_{i\geq
  r}\left\vert\partial^{\beta
  }\widehat{u}_{i}\left(\varphi,x\right)\right\vert,
~\forall~\left(\varphi,x\right)\in\A_{0}\times\overline{\omega}$$ 
which together with (\ref{complet-1.12.5}) and (\ref{complet-1.12.9}) implies
\begin{eqnarray}\label{complet-1.12.10}
\left\Vert\partial^{\beta}\widehat{U}_{r}\left(\varphi,\cdot\right)
\right\Vert&\leq&\sum\limits_{i\geq
  r}\left\Vert\partial^{\beta}\widehat{u}_{i}\left(\varphi,\cdot\right)\right\Vert\nonumber\\
&\leq&\sum\limits_{i\geq
  r}\left\Vert\partial^{\beta}\widehat{R}_{i}\left(\varphi,\cdot\right)\right\Vert\nonumber\\
&\leq&  2\sum\limits_{i\geq
  r}\widehat{\alpha}_{i+1}^{\bullet}\left(\varphi\right), ~\forall
~\varphi\in\A_{0}.
\end{eqnarray}
Now, define $\eta=\eta\left(\varphi\right):=\left[2i\left(\varphi
\right)\right]^{-1}, ~\forall ~\varphi\in\A_{0}$ then, since
$\varepsilon<\eta\left(\varphi\right)$ if
and only if $1-i\left(\varphi\right)\varepsilon>\frac{1}{2}$, we get
$$2\sum\limits_{i\geq
  r}\widehat{\alpha}_{i+1}^{\bullet}\left(\varphi\right)=\frac{2\left[i\left(\varphi\right)\varepsilon\right]^{r+1}}{1-i\left(\varphi\right) \varepsilon}<4\left(i\left(\varphi\right)\right)^{r+1}\varepsilon^{r+1}\leq 4\left(i\left(\varphi \right)\right)^{r+1}\varepsilon^{-0}.$$
Thus, from (\ref{complet-1.12.10}) it follows that
$$\left\Vert\partial^{\beta}\widehat{U}_{r}\left(\varphi_{\varepsilon},\cdot\right)\right\Vert\leq
c\left(\varphi\right)\varepsilon^{-0}, ~\forall ~\varphi\in\A_{0},~
\forall~ \varepsilon\in I_{\eta\left(\varphi\right)}~ 
\mbox{where} ~c\left(\varphi\right) :=4\left(i\left(\varphi\right)
\right)^{r+1}, ~\forall~\varphi \in \A_{0}$$
which obviously implies (\ref{complet-1.12.6}) since the finite sum
$(\widehat{f}_{\nu_{1}}+\sum\limits_{i\geq r}\widehat{u}_{i})$ is
moderated. Hence $U=\widehat{U}_{1}$ is moderated and
(\ref{complet-1.12.4}) is proved.
Finally, we will show that $f_{\nu _{i}}\underset{i\rightarrow\infty}{\longrightarrow} f$ in the topology $\mathcal{T}_{\overline{\omega}},_{b}$.
Fix an arbitrary $0$-neighborhood
$W_{\lambda,p}~\left(\lambda\in\Nset^{m},~
  p\in\Nset^{\ast}\right)$ in 
$\mathcal{G}\left(\overline{\omega}\right)$, we must show that there
exists $\theta\in\Nset$ verifying
\begin{equation}\label{complet-1.12.11}
t\geq \theta\Rightarrow (f-f_{\nu _{t+1}})\in W_{\lambda,p}.
\end{equation}
In view of (\ref{complet-1.12.2}) and (\ref{complet-1.12.3}) we have
\begin{equation}\label{complet-1.12.12}
\partial^{\varkappa}\widehat{f}\left(\varphi,x\right)=\partial^{\varkappa}\widehat{f}_{\nu _{1}}\left(\varphi,x\right)+\sum\limits_{i\geq 1}\partial^{\varkappa}\widehat{u}_{i}\left(\varphi,x\right)~ 
\forall
~\left(\varphi,x\right)\in\A_{0}\times\overline{\omega},~\forall~\varkappa\in\Nset^{m}. 
\end{equation}
For an arbitrary $t\in\Nset^{\ast }$, if $\varphi\in\A_{\nu
  _{t}}$ verifies $i\left(\varphi\right)\in \textsc{I}_{\nu
  _{t}^{-1}}$ [and therefore
$\left(\varphi,i\left(\varphi\right)\right)\in\A_{\nu
  _{i}}\times\textsc{I}_{\nu _{i}^{-1}}, ~\forall~ i=1,2,...,t$] it is clear
that
$$\partial^{\varkappa}\widehat{f}_{\nu _{1}}\left(\varphi,x\right)
+\sum\limits_{i=1}^{t}\partial^{\varkappa}\widehat{u}_{i}\left(\varphi,x\right)
=\partial^{\varkappa}\widehat{f}_{\nu_{t+1}}\left(\varphi,x\right), ~\left(x\in \overline{\omega },~\varkappa\in\Nset^{m}\right),$$
which shows that, for
$\left(\varphi,i\left(\varphi\right)\right)\in\A_{\nu_{t}}\times\textsc{I}_{\nu_{t}^{-1}}$, we can write (\ref{complet-1.12.12}) as
$$\partial^{\varkappa}\widehat{f}\left(\varphi,x\right)=\partial^{\varkappa}\widehat{f}_{\nu_{t+1}}\left(\varphi,x\right)+\sum\limits_{i\geq
  t+1}\partial^{\varkappa}\widehat{u}_{i}\left(\varphi,x\right) ~\left(x\in\overline{\omega},~\varkappa\in\Nset^{m}\right)$$ 
hence
$$\partial^{\varkappa}\left(\widehat{f}-\widehat{f}_{\nu _{t+1}}\right)
\left(\varphi,x\right)=\sum\limits_{i\geq t+1}\partial^{\varkappa}
\widehat{u}_{i}\left(\varphi,x\right),
~\forall~\left(\varphi,i\left(\varphi\right)\right)\in\A_{\nu_{t}}\times 
\textsc{I}_{\nu_{t}^{-1}},~\forall ~x\in\overline{\omega},~\forall~
\varkappa\in\Nset^{m}$$ and consequently 
$$\left\Vert
  \partial^{\varkappa}\left(\widehat{f}-\widehat{f}_{\nu_{t+1}}\right)\left(\varphi,\cdot\right)\right\Vert\leq\sum\limits_{i\geq t+1}\left\Vert 
\partial^{\varkappa}\widehat{u}_{i}\left(\varphi,\cdot\right)\right\Vert,
~\forall
~\left(\varphi,i\left(\varphi\right)\right)\in\A_{\nu_{t}}\times\textsc{I}_{\nu
  _{t}^{-1}},  ~\forall ~t\geq 1, ~\forall~\varkappa\in\Nset^{m},$$ 
which by (\ref{complet-1.12.5}) implies
\begin{equation}\label{complet-1.12.13}
\left\Vert\partial^{\varkappa}\left(\widehat{f}-\widehat{f}_{\nu_{t+1}}\right) 
\left(\varphi,\cdot\right)\right\Vert\leq\sum\limits_{i\geq
t+1}\left\Vert\partial^{\varkappa}\widehat{R}_{i}\left(\varphi,\cdot\right)\right\Vert,~ 
\forall~\left(\varphi,i\left(\varphi\right)\right)\in\A_{\nu_{t}}\times
\text{\textsc{I}}_{\nu _{t}^{-1}},t\geq
1,\varkappa\in\Nset^{m}.
\end{equation}
Now, we choose $\theta\in\Nset^{\ast}$ (see
(\ref{complet-1.12.11})) such that
\begin{equation}\label{complet-1.12.14}
\theta >p  \qquad \mbox{and} \qquad \lambda\leq\sigma_{\theta+1}
\left(\leq\sigma_{t+1}, ~\forall~ t\geq\theta\right).
\end{equation}
Clearly (\ref{complet-1.12.12}) and (\ref{complet-1.12.13}) hold for
any $\varkappa\in\Nset^{m}$ but now, in order to
prove (\ref{complet-1.12.11}), it suffices  to consider 
$\varkappa\leq\lambda.$ Next, we apply (\ref{complet-1.16.3}) (see
Lemma \ref{complet-1.16}) for an arbitrary fixed
$\varkappa\leq\lambda$, hence (by (\ref{complet-1.12.14})) we have
$\varkappa\leq\lambda\leq\sigma_{\theta+1}\leq\sigma_{t+1}, ~\forall
~t\geq\theta$, which shows that the term $\sigma_{i+1}$ in 
(\ref{complet-1.16.3}) should be replaced by $\sigma_{t+1}$ 
for all $t\geq\theta$. Therefore we can write
$$\left\Vert\partial^{\varkappa}\widehat{R}_{t}\left(\varphi,\cdot\right)\right\Vert\leq
2\widehat{\alpha}_{t+1}^{\bullet}\left(\varphi\right), ~\forall
~t\geq\theta, ~\forall~\varphi\in\A_{0}, ~\forall~ \varkappa\leq\lambda.$$
[Note that in (\ref{complet-1.16.3}) all the parameters vary freely: 
``$i\geq 1,~\varkappa\leq\sigma_{i+1},~\varphi\in\A_{0}$'' but in the  
above application of (\ref{complet-1.16.3}), $\varkappa$ is fixed by
the condition
$\varkappa\leq\lambda\leq\sigma_{\theta+1}$ (see (\ref{complet-1.12.14})) which
implies $\varkappa\leq\sigma_{t+1}, ~\forall~t\geq\theta$, hence the
above inequality is true for all $t\geq\theta$.]
and therefore (writing $\varphi_{\epsilon }$ instead of $\varphi$):
$$\left\Vert\partial^{\varkappa}\widehat{R}_{t}\left(
\varphi_{\epsilon},\cdot\right)\right\Vert\leq 2\left[i\left(
\varphi\right)\varepsilon\right]^{t+1}, ~\forall~ t\geq
\theta, ~\forall ~\varphi\in\A_{0},~\forall ~\varkappa\leq\lambda.$$
Hence, we can find an upper bound for the second member of
(\ref{complet-1.12.13}) in the following way: for 
$\varphi\in\A_{0},\varepsilon\in\textsc{I}_{\eta\left(\varphi\right)}~
(\eta\left(\varphi\right)$ was already defined by
$\eta\left(\varphi\right):=\left[2i\left(\varphi\right)\right]^{-1},
~\forall ~\varphi\in\A_{0}$) and $t\geq\theta$, we get
\begin{eqnarray}
\sum\limits_{i\geq t+1}\left\Vert\partial^{\varkappa}\widehat{R}_{i}
\left(\varphi_{\varepsilon},\cdot\right)\right\Vert&\leq&
\sum\limits_{i\geq\theta+1}\left\Vert\partial^{\varkappa}\widehat{R}_{i}\left(\varphi_{\epsilon},\cdot\right)\right\Vert\nonumber\\
&\leq&2\sum\limits_{i\geq
  \theta+1}\left[i\left(\varphi\right)\varepsilon\right]^{i+1}\nonumber\\
&=&\frac{2\left[i\left(\varphi\right)\varepsilon\right]}{1-i\left(\varphi\right)\varepsilon}^{\theta+2}\nonumber\\
&<& 4\left[i\left(\varphi\right)\varepsilon\right]^{\theta+2}\nonumber
\end{eqnarray}
and therefore we can write
\begin{equation}\label{complet-1.12.15}
\sum\limits_{i\geq t+1}\left\Vert\partial^{\varkappa
  }\widehat{R}_{i}\left(\varphi_{\varepsilon},\cdot\right)\right\Vert\leq
4\left[i\left(\varphi\right)\varepsilon \right]^{\theta +2}, ~\forall ~
\varphi\in\A_{0},~\forall~\varepsilon\in\text{\textsc{I}}_{\eta\left(\varphi\right)},~\forall ~t\geq \theta, ~\forall ~\varkappa\leq\lambda.
\end{equation}
Now, define
$\eta_{t}\left(\varphi\right):=\min\left(\eta\left(\varphi\right),\left[i\left(\varphi\right)\nu_{t}\right]^{-1}\right),
~\forall ~\varphi\in\A_{0}$, and $t\geq 1$. Then, by 
(\ref{complet-1.12.13}) and (\ref{complet-1.12.15}) we obtain:
\begin{equation}\label{complet-1.12.16}
\left\Vert\partial^{\varkappa}\left(\widehat{f}-\widehat{f}_{\nu_{t+1}}\right)\left(\varphi_{\varepsilon},\cdot\right)\right\Vert\leq
4\left[i\left(\varphi\right)\varepsilon\right]^{\theta+2}, ~\forall 
~t\geq\theta, ~\forall ~\varphi\in\A_{\nu_{t}},~\forall~\varepsilon 
\in \text{\textsc{I}}_{\eta_{t}\left(\varphi\right) } ~ \mbox{and}~\varkappa\leq\lambda.
\end{equation}
Next, note that since $\theta >p$ (see (\ref{complet-1.12.14})) and $\varepsilon \in \textsc{I}_{\eta _{t}\left(\varphi\right)}$ \Big(which implies
$i\left(\varphi\right)\varepsilon <\frac{1}{2}$\Big) it follows at
once that 
$4\left[i\left(\varphi\right)\varepsilon\right]^{\theta+2}<\left[i\left(\varphi\right)\varepsilon\right]^{p}$,
which by (\ref{complet-1.12.16}) implies
\begin{equation}\label{complet-1.12.17}
\left\Vert\partial^{\varkappa}\left(\widehat{f}-\widehat{f}_{\nu
_{t+1}}\right)\left(\varphi_{\varepsilon },\cdot\right)\right\Vert\leq
\left[i\left(\varphi\right)\varepsilon\right]^{p}, ~\forall ~ t\geq
\theta,~\forall ~\varphi\in\A_{\nu _{t}},
\forall~\varepsilon\in\text{\textsc{I}}_{\eta
  _{t}\left(\varphi\right)} ~\mbox{and} ~\varkappa\leq\lambda.
\end{equation}
Since the above inequality means that $\left[\widehat{\alpha}_{p}^{\bullet
}-\left(\widehat{f}-\widehat{f}_{\nu_{t+1}}\right)_{\varkappa}\right]
\left(\varphi_{\varepsilon}\right)\geq 0,$ the
statement (\ref{complet-1.12.17}) shows that for every 
$t\geq\theta$ we have proved that
{\tmsamp{
$$\left\vert 
\begin{array}{ll}
\exists ~N:=\nu_{t}\in\Nset~\mbox{such that}~\forall ~b>0
~\mbox{and}~ \forall ~\varphi\in\A_{N}=\A_{\nu_{t}}~\exists \eta:=\eta_{t}\left( \varphi\right)\in\text{\textsc{I}} \\ 
\mbox{verifying}~\left[\widehat{\alpha}_{p}^{\bullet}-\left(\widehat{f}-\widehat{f}_{\nu_{t+1}}\right)_{\varkappa}\right]\left(\varphi_{\varepsilon}\right)\geq
0>-\varepsilon^{b}, ~\forall ~\varepsilon \in\text{\textsc{I}}_{\eta}
~ \mbox{and}~ \varkappa\leq\lambda\end{array}\right.$$}}
and therefore (see [\cite{A-F-J}, Lemma 2.1 (ii)]),
(\ref{complet-1.12.11}) is proved. 
\end{proof}

\begin{Pro}\label{complet-1.17} {\tmsamp{The set
$$\mathcal{G}_{c}\left(\Omega\right):=\{f\in\mathcal{G}\left(\Omega\right)|\supp\left(f\right)\subset\subset\Omega\}$$
is a dense ideal of $\left(\mathcal{G}\left(\Omega \right),\mathcal{T}_{\Omega }\right)$.}}
\end{Pro}
\begin{proof} Let $\left(\chi_{l}\right)_{l\geq 1}$ be a
regularizing family associated to the fixed exhaustive
sequence $\left(\Omega _{l}\right)_{l\geq 1}$ for $\Omega$, that is,
$\chi_{l}\in\mathcal{D}\left(\Omega_{l+1}\right)$ and 
$\chi_{l}|_{\Omega_{l}}\equiv 1$ for each
$l\geq 1$. Then it will suffice  to show that for an arbitrary 
$f\in\mathcal{G}\left(\Omega\right)$ we have
$$\chi_{l}f\underset{l\rightarrow\infty}{\longrightarrow} f,$$
in the topology $\mathcal{T}_{\Omega}$. To this end fix an arbitrary
$0$-neighborhood $W_{\nu _{,}r}^{\beta }$ in
$\mathcal{G}\left(\Omega\right)$ ($\beta \in\Nset^{m}$ and $\nu,r\in\Nset^{\ast }$) and fix any representative $\widehat{f}$ of $f$.
Let $l_{0}\in\Nset^{\ast}$ be such that $\Omega
_{\nu}\subset\Omega_{l_{0}}\subset\Omega_{l}$ for all $l\geq l_{0}$, 
then $\left(\chi_{l}f\right)|_{\Omega _{l}}=f_{\Omega _{l}}, ~\forall
~l\geq l_{0}$ hence
(by setting $\widetilde{\chi}_{l}\left(\varphi,\cdot\right)
:=\chi_{l},~\forall ~l\geq l_{0}, ~\mbox{and}~\forall~\varphi\in\A_{0}$)
we have
$$\partial^{\sigma}\left(\widetilde{\chi}_{l}\widehat{f}-\widehat{f}\right)\left(\varphi,\cdot\right)|_{\Omega
  _{l}}\equiv 0, ~\forall ~l\geq l_{0}, ~\forall ~\sigma\in\Nset^{m}
~\mbox{and}~ \forall ~\varphi\in\A_{0}$$
which implies
$$\left\Vert\partial^{\sigma}\left(\widetilde{\chi}_{l}\widehat{f}-\widehat{f}\right)\left(\varphi,\cdot\right)\right\Vert
_{\nu}\equiv 0, ~\forall l\geq l_{0},~\forall~\sigma\leq\beta ~\mbox{and}
~\forall ~\varphi\in\A_{0}$$ and therefore $(\chi_{l}f-f)\in W_{\nu,r}^{\beta}, ~\forall ~l\geq l_{0}$.
\end{proof}

\begin{Lem}
\label{complet-1.18}
\begin{enumerate}
\item[$(a)$] {\tmsamp{The set $V_{r}[x_{0}]$ }}  (see [\cite{A-F-J}, 
Definition 2.7]) {\tmsamp{is bounded in $\left(\overline{\Kset},
\mathcal{T}\right)$ for
each $x_{0}\in\overline{\Kset}$ and each $r\in\Rset$.
\item[$(b)$] For given $\mu\in\overline{\Kset}$ and
  $V_{s}=V_{s}\left[0\right], \left(s>0\right)$ there is 
$N\in\Nset$ such that $\mu V_{s}\subset V_{-N}$.
\item[$(c)$] The set $\left\{V_{-N}|N\in\Nset\right\}$ is a fundamental system
  of bounded sets in $\left(\overline{\Kset},\mathcal{T}\right)$.
\item[$(d)$] If $X\subset\mathcal{G}\left(\overline{\omega}\right)$ 
satisfies the condition:
\begin{enumerate}
\item[$(B)$] $\forall~\beta\in\Nset^{m},~\exists ~N\in\Nset$ such that 
$\left\Vert u\right\Vert_{\mathcal{\sigma}}\leq
\alpha_{-N}^{\bullet}, ~\forall ~u\in X,~\forall ~\sigma \leq
\beta$. \newline 
Then $X$ is a bounded subset of 
$\left(\mathcal{G}\left(\overline{\omega}\right),
\mathcal{T}_{\overline{\omega}},_{b}\right)$.
\end{enumerate}}}
\end{enumerate}
\end{Lem}

\begin{proof}
In the proof of this result we will need the following
characterization of the relation $x\geq 0$ where
$x\in\overline{\Rset}$ and $\widehat{x}$ is any representative of
$x$ ( which is obviously equivalent to the conditions in [\cite{A-F-J}, Lemma 2.1])
{\tmsamp{
\begin{equation}\label{complet-1.18.1}
\left\vert 
\begin{array}{ll}
\exists ~N\in\Nset ~\mbox{such that}~\forall ~b_{0}>0,~\forall~
b>b_{0}~\mbox{and}~\forall ~\varphi\in \A_{N} \\ 
\exists ~\eta=\eta\left(b,\varphi\right)\in \text{\textsc{I}}
~\mbox{verifying} ~\widehat{x}\left(\varphi_{\varepsilon}\right)\geq-
\varepsilon^{b} ~\forall ~\varepsilon \in \text{\textsc{I}}_{\eta}.
\end{array}\right. 
\end{equation}}}

$(a)$ \textbf{case 1}: $x_{0}=0$ and $r>0$. Fix an arbitrary $V_{s}$
with $s>0$; it suffices  to show that there exists $t>0$ such that
\begin{equation}\label{complet-1.18.2}
V_{t}V_{r}\subset V_{s}.
\end{equation}
Indeed, it suffices to take $t>\max\left(s-r,0\right)$ by applying
[\cite{A-F-J}, Lemma 2.1 (i)] and then the proof of (\ref{complet-1.18.2}) 
follows at once.

\textbf{case 2}: $x_{0}=0$ and $r\leq 0$. Fix an arbitrary $V_{s}$
with $s>0$; choose
$t>\max \left(s-r,0\right)$; then the proof of (\ref{complet-1.18.2}) 
follows from (\ref{complet-1.18.1}) with $b_{0}:=-r>0$ 
if $\ r<0$ (the case $r=0$ is trivial).

\textbf{case 3}: $x_{0}\in \overline{\Rset}$ and
$r\in\Rset$.
Fix an arbitrary $V_{s}$ with $s>0$; then it suffices  to show that
there is a $\mathcal{T}$-neighborhood $W$ of $0$ such that
\begin{equation}\label{complet-1.18.2.2}
W\cdot V_{r}\left[x_{0}\right]\subset V_{s}
\end{equation}

It is easy to see that (\ref{complet-1.18.2.2}) follows from the
continuity of the addition and multiplication in 
$\left(\overline{\Kset},\mathcal{T}\right) $ and from the fact
that $V_{r}$ is bounded for each $r\in\Rset$ ~({\bf{cases 1 and 2}}).

$(b)$ If $\widehat{\mu}$ is any representative of $\mu$ then there is $N\in 
\Nset$ such that $\forall ~\varphi\in \A_{N}~\exists ~c>0$ and
$\eta\in \textsc{I}$ verifying 
$$\left\vert\widehat{\mu}\left(\varphi_{\varepsilon}\right)\right\vert
\leq c\varepsilon^{-N}, ~\forall~\varepsilon \in \textsc{I}_{\eta }.$$
Now, if $\widehat{x}$ is any representative of an arbitrary $x\in V_{s}$ it suffices to apply [\cite{A-F-J}, Lemma 2.1 (i)] to $\widehat{x}$ and to 
$\widehat{\mu }\widehat{x}$ for to get $\left\vert\mu
  x\right\vert\leq\alpha _{-N}^{\bullet }$. 

$(c)$ Let $B$ be a bounded set in $\left(\mathcal{G}\left(\overline{\omega 
}\right),\mathcal{T}_{\overline{\omega },{b}}\right)$ then, for a
given $V_{s}$ with $s>0$, there is $V_{t}$ with $t>0$ such that
$$V_{t}B\subset V_{s}.$$
Since $0\in\overline{\Inv\left(\overline{\Kset}\right)}$ there is 
$\lambda^{-1}\in V_{t}\cap \Inv\left(\overline{\Kset}\right)$ and
then, from the above inclusion, it follows that 
$\lambda^{-1}B\subset V_{s}$ and therefore $B\subset\lambda V_{s}$ and the
conclusion follows from $(b)$.

$(d)$ Fix an arbitrary $W_{\beta,t}~(\beta\in\Nset^{m}, ~t\in\Rset_{+}^{\ast })$, then we must show that there exists $r>0$ such that
\begin{equation}\label{complet-1.18.3}
V_{r}X\subset W_{\beta,t}
\end{equation}
From the assumption $(B)$, for the above fixed
$\beta\in\Nset^{m}, ~\exists~ N\in\Nset$ such that
$$\left\Vert u\right\Vert_{\sigma}\leq \alpha_{-N}^{\bullet}, ~
\forall~ u\in X, ~\forall ~\sigma\leq\beta.$$
We set $r:=N+t$. Fix $\lambda\in V_{r}, ~u\in X$ and representatives
$\widehat{\lambda}$ and $\widehat{u}$ of $\lambda$ and $u$
respectively. From (\ref{complet-1.18.1}), $\lambda \in V_{r}$ means that
{\tmsamp{
\begin{equation}\label{complet-1.18.4}
\left\vert\begin{array}{ll}
\exists ~N_{1}\in\Nset ~\mbox{such that}~\forall ~b_{0}^{\prime }>0,~\forall ~b^{\prime}>b_{0}^{\prime }~\mbox{and}~ 
\forall ~\varphi\in\A_{N_{1}}~\exists ~\eta_{1}(2b^{\prime
},\varphi)\in \text{\textsc{I}}\\
\mbox{verifying}~ 
\left\vert\widehat{\lambda}\left(\varphi_{\varepsilon}\right)
\right\vert\leq i\left(\varphi\right)^{r}\varepsilon^{r}+\varepsilon
^{2b^{\prime }}~\forall ~\varepsilon\in\text{\textsc{I}}_{\eta
  _{1}}.\end{array}\right.  
\end{equation}}}
On the other hand, from the condition $(B)$, the relation $u\in X$ means that
{\tmsamp{
\begin{equation}\label{complet-1.18.5}
\left\vert \begin{array}{ll}
\exists ~N^{\prime}\in\Nset~\mbox{such that}~\forall
~b_{0}^{\prime\prime}>0~ \forall
~b^{\prime\prime}>b_{0}^{\prime\prime}
~\mbox{and} ~
\forall ~\varphi\in \A_{N^{\text{ }^{\prime }\text{ }}}\exists 
~\eta =\eta\left(2b^{\prime \prime},\varphi\right)\in
\text{\textsc{I}}\\
\mbox{verifying} ~ 
\left\Vert\partial^{\sigma}\widehat{u}\left(\varphi_{\varepsilon},\cdot\right) 
\right\Vert\leq
i\left(\varphi\right)^{-N}\varepsilon^{-N}+\varepsilon^{2b^{\prime
    \prime }} ~\forall ~\varepsilon \in \text{\textsc{I}}_{\eta
},~\forall ~\sigma\leq \beta.\end{array}\right.
\end{equation}}}
Next we are going to prove (\ref{complet-1.18.3}), that is, 
$\lambda u\in W_{\beta ,t}.$ We define
$N_{0}:=\max\left(N^{\prime},~N_{1}\right)$. 
Fix $b>b_{0}:=\max\left(r,t\right),~\varphi\in \A_{0}$ and
set $$\eta_{\ast}=\eta_{\ast}\left(b,\varphi\right):=\min\left(\eta\left(2b,\varphi\right),\eta_{1}\left(2b,\varphi\right)\right)\in \textsc{I}.$$
We apply (\ref{complet-1.18.4}) with $b_{0}^{\prime }:=\frac{1}{2}\left(N+b\right)$ and (\ref{complet-1.18.5}) with  
$b_{0}^{\prime \prime}:=\frac{1}{2}\left(b-r\right)$. From the definition
of $b$ it is clear that $b>b_{0}^{\prime }$ and
$b>b_{0}^{\prime \prime }$ hence $\eta_{\ast}$ is well defined. Now it is
trivial to see that there exists
$\eta\in\textsc{I},~\eta\leq\eta_{\ast}$ verifying
$$\left\Vert\partial^{\sigma}\left(\widehat{\lambda}\widehat{u}\right)\left(\varphi_{\varepsilon},\cdot\right)\right\Vert\leq
i\left(\varphi\right)^{t}\varepsilon^{t}+\varepsilon^{b}, ~\forall ~\varepsilon \in \textsc{I}_{\eta},~~\forall~\sigma\leq\beta,$$
that is, $\lambda u\in W_{\beta,t}$.
\end{proof}
\section{Two auxiliary results}\label{auxiliary}
\hspace{4.1mm} Fix a bounded open subset $\Omega$ of
$\Rset^{m},~T\in\Rset_{+}^{\ast}$ 
and set $Q:=\Omega \times \left] 0,T\right[\subset \Rset^{m+1}.$
For given $f\in\mathcal{G}\left(\overline{Q}\right)$ and $t_{0}\in\left[0,T\right]$ we must give a sense to the
"restriction"
$$R:=f|_{\overline{\Omega}\times\left\{t_{0}\right\}}$$
showing that $R$ can be identified naturally to an element of $\mathcal{G}\left(\overline{\Omega}\right)$.

Note that $R$ does not a priori make sense since\ $int\left(\overline{\Omega}\times\left\{t_{0}\right\}\right)=\varnothing$ and hence
$\overline{\Omega}\times\left\{t_{0}\right\}$ is not a quasi regular set (see [\cite{A}, Definition 1.1]).

Fix a representative $\widehat{f}\in \mathcal{E}_{M}\left[\overline{Q}\right]$ of $f$; then the restriction

$$\widehat{f}|_{\A_{0}\left(m+1\right)\times\overline{\Omega}\times\left\{t_{0}\right\}}$$
is well defined since 
$\A_{0}\left(m+1\right)\times\overline{\Omega}\times\left\{t_{0}\right\}\subset \A_{0}\left(m+1\right)\times\overline{Q}$. Therefore,
we can define (see the definition of \textsc{I}$_{m+1}^{m}$\ in
[\cite{ab}, Proposition 1.7 (e)]) the function
$$f_{t_{0}}^{\ast}:\left(\psi,x\right)\in
\A_{0}\left(m\right)\times\overline{\Omega}\mapsto\widehat{f}\left(\text{\textsc{I}}_{m+1}^{m}\left(\psi\right),x,t_{0}\right) \in 
\Kset.$$

With the above notation we have
\begin{Lem}\label{aux-2.1}
\begin{enumerate}
\item[$(a)$] {\tmsamp{$f_{t_{0}}^{\ast}\in\mathcal{E}_{M}\left[\overline{\Omega}\right]$; 
\item[$(b)$] If $\widehat{g}$ is another representative of $f$
and $g_{t_{0}}^{\ast}$ is the function defined similarly as $f_{t_{0}}^{\ast}$ by
changing $\widehat{f}$ by $\widehat{g}$, then  $f_{t_{0}}^{\ast}-g_{t_{0}}^{\ast}$
$\in \mathcal{N}\left[\overline{\Omega}\right]$. Therefore, if $f_{t_{0}}:=\cl\left(f_{t_{0}}^{\ast}\right)\in\mathcal{G}\left(\overline{\Omega}\right)$, then we have a $\Kset$-linear map
$$r_{\overline{\Omega}}^{\overline{Q}}:f\in\mathcal{G}\left(\overline{Q}\right)\mapsto f_{t_{0}}=\cl\left(f_{t_{0}}^{\ast}\right)\in\mathcal{G}\left(\overline{\Omega}\right).$$
\item[$(c)$] The map $r_{\overline{\Omega}}^{\overline{Q}}$ is
  $(\mathcal{T}_{\overline{Q},b}-\mathcal{T}_{\overline{\Omega},b})$--continuous.}}
\end{enumerate}
\end{Lem}
\begin{proof} $(a)$ Fix $K\subset\subset\overline{\Omega}$
and $\varkappa\in\Nset^{m}$ then $\varkappa^{\prime}:=\left(\varkappa,0\right)\in\Nset^{m}\times\Nset$ and $K\times\left\{t_{0}\right\}\subset\subset\overline{Q}$. For the fixed
representatives $\widehat{f}$ and $\ f_{t_{0}}^{\ast}$ of $f$ and $f_{t_{0}}$ respectively, 
\begin{equation}\label{aux-2.1.1}
\partial^{\varkappa^{\prime}}\widehat{f}\left(\varphi,x,t_{0}\right)=\partial^{\varkappa}f_{t_{0}}^{\ast}\left(\text{\textsc{I}}_{m}^{m+1}\left(\varphi\right),x\right),~\forall~x\in\overline{\Omega},~\forall~\varphi\in \A_{0}\left(m+1\right),
\end{equation}
which implies that the moderation of $f_{t_{0}}^{\ast }$ follows at once
from the moderation of $\widehat{f}$.

$(b)$ Follows immediately from (\ref{aux-2.1.1}).

$(c)$ Fix any $\mathcal{T}_{\overline{\Omega},b}$-neighborhood $W_{\sigma,r}\left( \sigma\in\Nset^{m},~r>0\right)$ of $0$ in $\mathcal{G}\left(\overline{\Omega}\right)$.
We set
$$\sigma^{\prime}:=\left(\sigma,0\right)\in\Nset^{m}\times\Nset \ \mbox{and}\  ~s>\frac{r\left(m+1\right)}{m}$$
and we shall show that
\begin{equation}\label{aux-2.1.2}
f\in W_{\sigma^{\prime },s}^{\ast}\Rightarrow f_{t_{0}}\in W_{\sigma,r}
\end{equation}
where $W_{\sigma^{\prime},s}^{\ast}$ denote the subset
$W_{\sigma^{\prime },s}$ of
$\mathcal{G}\left(\overline{Q}\right)$. Since in this proof we shall
work with two different dimensions we shall need the elements 
$\alpha_{r}^{\bullet}\in\overline{\Kset}\left(m\right)$
represented, as usual, by
$$\widehat{\alpha}_{r}^{\bullet}:\psi\in \A_{0}\left(m\right)
\mapsto i\left(\psi\right)^{r}\in\Rset _{+}^{\ast}$$
and also the elements $\beta_{s}^{\bullet}\in\overline{\Kset}\left(m+1\right)$ represented by
$$\widehat{\beta}_{s}^{\bullet}:\varphi\in
\A_{0}\left(m+1\right)\mapsto i\left(\varphi\right)^{s}
\in\Rset_{+}^{\ast}.$$
The relation $f\in W_{\sigma^{\prime},s}^{\ast}$ means that
\begin{equation}\label{aux-2.1.3}
\left\Vert f\right\Vert_{\varkappa^{\prime }}\leq\beta_{s}^{\bullet},
~\forall~ \varkappa^{\prime}\leq\sigma^{\prime}.
\end{equation}
Note that $$\varkappa^{\prime}\leq\sigma^{\prime}=\left(\sigma
,0\right)\Rightarrow \varkappa^{\prime}=\left(\varkappa_{1},\dots\varkappa_{m},0\right)=\left(\varkappa,0\right)$$ where
$\varkappa:=\left(\varkappa_{1},\dots,\varkappa_{m}\right)$ hence 
$\varkappa^{\prime}\leq\sigma^{\prime}\Longleftrightarrow \varkappa 
\leq\sigma$, therefore we can write (\ref{aux-2.1.3}) in the following way
\begin{equation}\label{aux-2.1.3.3}
\left\Vert f \right\Vert_{\varkappa^{\prime }}\leq\beta_{s}^{\bullet},
~\forall~\varkappa\leq\sigma . 
\end{equation}

From the definition of $f_{t_{0}}^{\ast}$ we get
\begin{equation}\label{aux-2.1.4}
\left\Vert\partial^{\varkappa}f_{t_{0}}^{\ast}\left(\psi,-\right)\right\Vert_{\overline{\Omega}}\leq\left\Vert\partial^{\varkappa
    ^{\prime
    }}\widehat{f}\left(\varphi,-,-\right)\right\Vert_{\overline{Q}},
~\forall ~\varphi\in \A_{0}\left(m+1\right),
~\psi:=\textsc{I}_{m}^{m+1}\left(\varphi\right),~\varkappa^{\prime}\leq
\sigma^{\prime }.
\end{equation}
Now, it is easily seen that the function (see [\cite{ab}, Lemma 3.1.1]):
$$\psi\in \A_{0}\left(m\right)\mapsto\left\Vert\partial^{\varkappa^{\prime}}\widehat{f}\left(\text{{\small I}}_{m+1}^{m}\left(\psi \right),-,-\right)\right\Vert _{\overline{Q}}\in\Rset_{+}$$
is a representative of $\left\Vert f_{{t_0}}\right\Vert_{\varkappa^{\prime}}$ and since the function
$$\psi \in \A_{0}\left(m\right)\mapsto\left\Vert \partial^{\varkappa}f_{t_{0}}^{\ast}\left(\psi ,-\right)\right\Vert_{\overline{\Omega}}\in\Rset_{+}$$
is a representative of $\left\Vert f_{t_{0}}\right\Vert _{\varkappa }$, it
is clear from (\ref{aux-2.1.4}) and [\cite{A-F-J}, Lemma 2.1(iii)] that
\begin{equation}\label{aux-2.1.5}
\left\Vert f_{t_{0}}\right\Vert_{\varkappa}\leq\left\Vert f
\right\Vert_{\varkappa ^{\prime }}, ~\forall ~\varkappa^{\prime}\leq
\sigma^{\prime }\left(\Longleftrightarrow\varkappa\leq\sigma\right).
\end{equation}
Next, note that
$\widehat{\beta}_{s}^{\bullet}\in\mathcal{E}_{M}\left(\Kset,m+1\right)$
hence
$^{\ast}\beta_{s}^{\bullet}:=\widehat{\beta}_{s}^{\bullet}\circ\textsc{I}_{m+1}^{m}\in\mathcal{E}_{M}\left(\Kset,m\right)$ 
is a representative of $\beta _{s}^{\bullet }$ in $\mathcal{E}_{M}\left( 
\Kset,m\right)$ [since $^{\ast}\beta_{s}^{\bullet}$ is a
representative of
$B_{m}^{m+1}\left(\beta_{s}^{\bullet}\right)\in\overline{\Kset}
\left(m\right)$, where $B_{m}^{m+1}:\overline{\Kset}\left(m+1\right) 
\rightarrow\overline{\Kset}\left(m\right)$ is the natural
isomorphism induced by the map 
$$\widetilde{B}_{m}^{m+1}:w\in\mathcal{E}_{M}\left(\Kset,m+1\right)\mapsto 
w\circ\textsc{I}_{m+1}^{m}\in \mathcal{E}_{M}\left(\Kset,m\right),$$ see [\cite{ab}, Lemma 3.1.1]].
 Hence, in order to prove that
$\beta _{s}^{\bullet}\leq\alpha_{r}^{\bullet}$ we can work with the
representatives $\widehat{\alpha}_{r}^{\bullet}$ and
$$^{\ast}\beta _{s}^{\bullet}:\psi\in \A_{0}\left(m\right)\mapsto i\left(\text{\textsc{I}}_{m+1}^{m},\left(\psi\right)\right)^{s}\in\Rset_{+},$$
by applying [\cite{A-F-J}, Lemma 2.1(ii)]. Indeed, take $N:=0$, fix $b>0$ and
$\psi \in \A_{0}\left(m\right)$ arbitrary. Since our hypothesis on $s$
implies $$\lambda:=\frac{sm}{m+1}-r>0$$ 
it is obvious that there exists $\eta =\eta\left(\psi\right)\in
\textsc{I}$ such that
$$L\left(\psi,\varepsilon\right):=i\left(\psi\right)^{r}-i\left(\text{\textsc{I}}_{m+1}^{m}\left(\psi\right)\right)^{s}\varepsilon^{\lambda}\geq
0, ~\forall ~\varepsilon \in \textsc{I}_{\eta}.$$
 Since
$\left(\widehat{\alpha}_{r}^{\bullet}~-~^{\ast}\beta_{s}^{\bullet}\right)
\left(\psi_{\varepsilon}\right)=\varepsilon^{r}L\left(\psi,\varepsilon\right),
\forall ~\varepsilon >0$, it follows that
$$\left(\widehat{\alpha}_{r}^{\bullet}~-~^{\ast}\beta_{s}^{\bullet}\right)\left(\psi_{\varepsilon}\right)\geq
0>-\varepsilon^{b}, ~\forall ~\varepsilon \in \textsc{I}_{\eta}$$
hence $\beta_{s}^{\bullet}\leq\alpha_{r}^{\bullet}$. This inequality
together (\ref{aux-2.1.3}) and (\ref{aux-2.1.5}) shows that
$\left\Vert f_{t_{0}}\right\Vert_{\varkappa}\leq\alpha_{r}^{\bullet },
~\forall ~\varkappa\leq\sigma$, that is, $f_{t_{0}}\in W_{\sigma,r}$ 
and (\ref{aux-2.1.2}) is proved.
\end{proof}

The remainder of this section is devoted to the definition of a
suitable topology on the algebra $\mathcal{G}\left(\partial\Omega\times
\left[0,T\right]\right)$. First, we present a 
definition of $\mathcal{G}\left(\partial\Omega\times X\right)$, where 
$X$ is a quasi-regular set in $\Rset^{n}$ and $\Omega$ is
an open subset of $\Rset^{m}$ (this is necessary since in 
[\cite{A}, Definition 3.7],  two conditions 
obviously needed in the definition were omitted).
The set $\mathcal{E}_{M}\left[\partial\Omega\times X\right]$ is
defined as the set of all functions (where $\A_{0}=\A_{0}\left(m+n\right)$)
$$\widehat{u}:\A_{0}\times \partial\Omega\times X\rightarrow\Kset$$
verifying the conditions below:
{\tmsamp{
\begin{enumerate}
\item[(M$_{\text{\textsc{I}}}$)] $\left(\xi\in\partial\Omega\mapsto 
\widehat{u}\left(\varphi,\xi,t\right)\in\Kset\right)\in\mathcal{C}
\left(\partial\Omega,\Kset\right), ~\forall ~\varphi\in
\A_{0},~\forall ~t\in X$
\item[(M$_{\text{\textsc{II}}}$)] $\left(t\in X\mapsto \widehat{u}\left(
\varphi,\xi,t\right)\in\Kset\right)\in\mathcal{C}^{\infty}\left(X;\Kset\right),
~\forall ~\varphi\in \A_{0},~\forall ~\xi\in\partial\Omega$
\item[(M$_{\text{\textsc{III}}}$)] $\left\vert 
\begin{array}{ll} 
\forall ~\varkappa
\in\Nset^{n},~\forall~K\subset\subset X~\mbox{and}~
H\subset\subset\partial\Omega~\exists ~N\in\Nset~\mbox{such that}\\
~\forall ~\varphi\in \A_{N} ~\exists ~c=c\left(\varphi\right)>0~\mbox{and}~\eta 
 =\eta\left(\varphi\right)\in \text{\textsc{I }} \\ 
\mbox{verifying}~\left\vert\partial_{t}^{\varkappa}\widehat{u}\left(
\varphi_{\varepsilon},\xi,t\right)\right\vert\leq c\varepsilon^{-N},
~\forall ~ t\in K,~\xi\in H~\mbox{and}~ \varepsilon\in 
\text{\textsc{I}}_{\eta}.\end{array}\right.$
\end{enumerate}}}

Clearly, the set $\mathcal{E}_{M}[\partial\Omega\times X]$ with the
pointwise operations is a $\overline{\Kset}$-algebra
and the set $\mathcal{N}\left[\partial\Omega\times X\right]$ of all 
$\widehat{u}\in\mathcal{E}_{M}\left[\partial\Omega\times X\right]$
verifying the condition:
{\tmsamp{
$$\mbox{(N)}\qquad\left\vert\begin{array}{ll}
\forall~\varkappa\in\Nset^{n},~\forall ~K\subset\subset X
~\mbox{and}~\forall ~ H\subset\subset\partial\Omega~\exists ~N\in\Nset
~\mbox{and}~ \gamma\in\Gamma\\ 
\mbox{such that} ~\forall ~q\geq N ~\mbox{and} ~\forall ~\varphi\in
\A_{q}~ \exists ~c=c\left(\varphi\right)>0 ~\mbox{and}~\eta=\eta\left(\varphi 
\right)\in \text{\textsc{I}} \\ 
~\mbox{verifying}~\left\vert\partial_{t}^{\varkappa}\widehat{u}\left(\varphi_{\varepsilon},\xi,t\right)\right\vert\leq
c\varepsilon^{\gamma\left(q\right)-N}, ~\forall ~t\in K, \xi\in H ~\mbox{and} 
~\varepsilon\in\text{\textsc{I}}_{\eta}\end{array}\right.$$}} 
is an ideal of $\mathcal{E}_{M}\left[\partial\Omega\times X\right]$.
This allows to set
$$\mathcal{G}\left(\partial\Omega\times
  X\right):=\frac{\mathcal{E}_{M}
\left[\partial\Omega\times X\right]}{\mathcal{N}\left[\partial\Omega\times 
 X\right]}.$$ 

In the sequel we restrict attention to the present  case: $\Omega$ is a bounded open
subset of $\Rset^{m}$ and
$X:=\left[0,T\right]\subset\Rset_{+}$ 
$\left(T>0\right)$. Therefore, since $\partial\Omega$ and $X$ are
compact sets, the conditions (M$_{\text{\textsc{III}}}$) and (N) can be
simplified for it is enough to work with $K=X$ and
$H=\partial\Omega$. In what follows, for the sake of simplicity, we set
$$Q:=\Omega\times\left]0,T\right[ ~\mbox{and}~
Q^{\ast}:=\partial\Omega\times\left[0,T\right].$$
and we will define a topology on $\mathcal{G(}Q^{\ast})$. For a given 
$u\in \mathcal{G}\left(Q^{\ast}\right)$ let $\widehat{u}$ be any
representative of $u$ then, for every $\nu\in\Nset$ the function
$$\widehat{u}_{\left(\nu\right)}:\varphi \in \A_{0}\mapsto \|\partial_{t}^{\nu }\widehat{u}\left(\varphi,\xi,t\right)\|_{Q^{\ast}}:=\underset{\left(\xi,t\right) \in Q^{\ast}}{\sup}\left\vert\partial_{t}^{\nu}\widehat{u}\left(\varphi,\xi,t\right)\right\vert\in\Rset_{+}$$
is moderate. If $u^{\ast}\in\mathcal{E}_{M}\left[Q^{\ast}\right]$ and 
$u_{\left(\nu\right)}^{\ast}$ is defined as $\widehat{u}_{\left(\nu
\right)}$ with $\widehat{u}$ replaced by $u^{\ast }$, it follows at
once (see [\cite{A-F-J}, Lemma 3.1]) that
$$(\widehat{u}-u^{\ast})\in \mathcal{N}\left[Q^{\ast}\right]\Rightarrow
\left(\widehat{u}_{\left(\nu\right)}-u_{\left(\nu\right)}^{\ast}\right)
\in\mathcal{N}\left(\Kset\right)$$ which shows that
$\cl\left(\widehat{u}_{\left(\nu\right)}\right)\in\overline{\Rset}_{+}$ 
is independent of the chosen representative chosen of $u$. Therefore, for each
$\nu\in\Nset$, we have a function
$$p_{\left(\nu\right)}:u\in\mathcal{G}\left(Q^{\ast}\right)\mapsto \cl
\left(\widehat{u}_{\left(\nu\right)}\right)\in\overline{\Rset}_{+}$$
where $\widehat{u}$ is any representative of $u$. Clearly $p_{\left(\nu\right)}$ is a $G$-seminorm
(see [\cite{A-F-J}, Definition A.1]), that is, for all $u,v\in\mathcal{G}\left(
Q^{\ast}\right)$ and $\lambda\in\overline{\Kset}$ we have
$$p_{\left(\nu\right)}\left(u+v\right)\leq
p_{\nu}\left(u\right)+p_{\left(\nu\right)}\left(v\right)
\qquad\mbox{and} \qquad p_{\left(\nu\right)}\left(\lambda u\right)=\left\vert\lambda\right\vert p_{\left(\nu\right)}\left(u\right).$$

It follows that the set $\mathcal{B}_{Q^{\ast}}$ of the sets
$\left(\nu\in\Nset,s\in\Rset\right)$:
$$N_{\nu,s}:=\left\{u\in\mathcal{G}\left(Q^{\ast}\right)|p_{\left(l\right)}\left(u\right)\leq\alpha_{s}^{\bullet}, ~\forall ~ l\leq \nu\right\}$$
is a fundamental system of $0$-neighborhoods for a topology 
$\mathcal{T}_{Q^{\ast}}$ on $\mathcal{G}\left(Q^{\ast}\right)$
which is compatible with the $\overline{\Kset}$-algebra structure
of $\mathcal{G}\left(Q^{\ast}\right)$ (the proof of
this statement, which consists in showing that $\mathcal{B}_{Q^{\ast }}$
satisfies the seven condition of [\cite{A-F-J}, Proposition 1.2~ 
($2^{\underline{o}}$)], is easy but rather tedious and we do not give it here).

Now, fix $w\in\mathcal{G}\left(\overline{Q}\right)$ and any
representative $\widehat{w}\in\mathcal{E}_{M}\left[\overline{Q}\right]$
of $w$. Then the restriction 
$$\widehat{u}:=\widehat{w}|_{\A_{0}\times Q^{\ast}}$$
is well defined and 
$\widehat{u}\in\mathcal{E}_{M}\left[ Q^{\ast}\right]$. Moreover, if 
$w^{\ast}\in\mathcal{E}_{M}\left[\overline{Q}\right] ~\mbox{and} ~ 
u^{\ast }:=w^{\ast}|_{\A_{0}\times Q^{\ast }}$, it is clear that
$$(\widehat{w}-w^{\ast})\in\mathcal{N}\left[\overline{Q}\right]\Rightarrow 
(\widehat{u}-u^{\ast})\in\mathcal{N}\left[Q^{\ast}\right]$$

which shows that we get a natural homomorphism of
$\overline{\Kset}$-algebras 
$$\rho=\rho_{Q}:w\in\mathcal{G}\left(\overline{Q}\right)\mapsto 
\cl\left(\widehat{w}|_{\A_{0}\times
    Q^{\ast}}\right)\in\mathcal{G}\left(Q^{\ast}\right),$$
where $\widehat{w}$ is any representative of $w$. 
In the result below we assume that
$\mathcal{G}\left(Q^{\ast }\right)$
(resp. $\mathcal{G}\left(\overline{Q}\right) )$ is endowed with the
topology $\mathcal{T}_{Q^{\ast }}$
(resp. $\mathcal{T}_{\overline{Q},b}$, see Theorem \ref{complet-1.7}).

\begin{Lem}\label{aux-2.2}{\tmsamp{
The above map $\rho=\rho_{Q}$ is continuous.}}
\end{Lem}
\begin{proof} Since the topology $\mathcal{T}_{\overline{Q},b}$ (resp. $
\mathcal{T}_{Q^{\ast}}$) is defined by the set (see Theorem \ref{complet-1.7})
$$\mathcal{B}_{\overline{Q},b}:=\{W_{\sigma,r}|\sigma=\left(\sigma^{\prime},
\nu \right)\in\Nset^{m}\times\Nset ~\mbox{and}
~r\in\Nset^{\ast}\}$$ (resp $\mathcal{B}_{Q^{\ast
  }}:=\{N_{\nu,s}|\nu\in\Nset ~\mbox{and}
~s\in\Nset^{\ast}\}$) 
it suffices to show that for an arbitrary given
$N_{\nu,s}\in\mathcal{B}_{Q^{\ast }}$ 
there is $W_{\sigma,r}\in\mathcal{B}_{\overline{Q},b}$ such that
\begin{equation}\label{aux-2.2.1}
w\in W_{\sigma,r}\Rightarrow \rho\left(w\right)\in N_{\nu,s}
\end{equation}
Fix $N_{\nu ,s}$ arbitrary in $\mathcal{B}_{Q^{\ast }}$ then by
defining $\sigma:=\left(0,\nu\right)\in\Nset^{m}\times\Nset
~\mbox{and}~r:=s$ it is easy to check that
$$w\in W_{\left(0,\nu\right),s}\Rightarrow\rho\left(w\right)\in N_{\nu,s}.$$
\end{proof}

\section{An initial boundary value problem}\label{ibvp}
\hspace{4.1mm} In this section we shall use the following notation:\\
{\tmsamp{
\begin{enumerate}
\item $\Rset[x]$ denotes the ring of the polynomials in one variable with real 
coefficients;
\item $\Omega$ is a non-void bounded open subset of $\Rset^m, 
T\in\Rset_{+}^{\ast }$ and hence $Q:=\Omega\times ]0,T[$ is a bounded
open subset of $\Rset^{m+1}$;
\item For the definition of
  $\mathcal{C}^{\infty}\left(\overline{Q}\right)$ 
see [\cite{A}, section1];
\item If $u\in\mathcal{C}^{\infty}\left(\overline{Q}\right)$
(resp. $u_{0}\in\mathcal{D}\left(\Omega\right)$) we set
$\left\Vert
  u\right\Vert_{\mathcal{C}^{k}\left(\overline{Q}\right)}:=\underset{\left\vert\sigma \right\vert\leq k}{\sum}\left\Vert\partial^{\sigma }u\right\Vert_{\overline{Q}}$ (resp. $\left\Vert u_{0}\right\Vert_{\mathcal{C}^{2k+1}\left(\overline{\Omega}\right)}:=\underset{\left\vert\tau\right\vert\leq 2k+1}{\sum}\left\Vert\partial^{\tau }u_{0}\right\Vert_{\overline{\Omega}})$;
\item IBVP is an abbreviation for "initial boundary value problem "
\end{enumerate}}}

Next, for the benefit of the reader, we shall begin by presenting four
results of [\cite{cl}] which we will need for to solve our problem.

\begin{Lem}([\cite{cl}, Lemma 1])\label{ibvp-3.1} {\tmsamp{For any $k\in \Nset$ there exists
  $P_{k}\in\Rset[x]$ such that, if $u_{0}\in \mathcal{D}\left( \Omega \right)$
then the solution $u\in \mathcal{C}^{\infty
}\left(\overline{Q}\right)$ of
\begin{equation}\label{ibvp-3.1.1}
\left\{ 
\begin{array}{ll}
u_{t}-\Delta u+u^{3}=0~\mbox{in}~ Q \\ 
u|_{\Omega \times \left\{ 0\right\}} =u_{0}\\
u|_{\partial \Omega \times \left] 0,T\right[} =0 \\ 
\end{array}\right.
\end{equation}
satisfies $\left\Vert \text{ }u\text{ }\right\Vert _{
\mathcal{C}^{k}\left( \overline{Q}\right)}\leq P_{k}\left(\left\Vert
u_{0}\right\Vert_{\mathcal{C}^{2k+1}\left(\overline{\Omega}\right)}\right)$.}}
\end{Lem}

\begin{The}([\cite{cl}, Theorem 1])\label{ibvp-3.2} {\tmsamp{For any
      given $u_{0}\in
  \mathcal{G}_{c}\left( \Omega,\Rset\right)$ there exists a solution 
$u\in \mathcal{G}\left(\overline{Q}\right)$ of 
\begin{equation}\label{ibvp-3.2.1}
\left\{\begin{array}{ll}
u_{t}-\Delta u+u^{3}=0~\mbox{in} ~\mathcal{G}\left(Q\right) \\ 
u|_{\overline{\Omega }\times \left\{0\right\}}
=u_{0}~\mbox{in}~\mathcal{G}
\left(\overline{\Omega}\right)\\ 
u|_{\partial \Omega \times \left[ 0,T\right]}
=0~\mbox{in}~\mathcal{G}\left(\partial\Omega\times\left[0,T\right]\right).
\end{array}\right.
\end{equation}}}
\end{The}

\begin{Lem}([\cite{cl}, Lemma 2])\label{ibvp-3.3} {\tmsamp{Let 
$v\in \mathcal{C}^{\infty}\left(\overline{Q}\right)$ be a solution of 
\begin{equation}\label{ibvp-3.3.1}
\left\{\begin{array}{ll}
v_{t}-\Delta v+a_{0}v=f ~\mbox{in} ~Q~(\mbox{with}~
f\in \mathcal{C}^{\infty}\left(\overline{Q}\right)) \\ 
v\left(x,0\right)=g\left(x\right)~\mbox{on}~
\Omega ~(\mbox{with}~g\in \mathcal{D}\left(\Omega\right))\\ 
v\left(x,t\right)=h\left(x,t\right)~\mbox{in}~\partial \Omega
\times \left[0,T\right]~(\mbox{with}~h\in \mathcal{C}^{\infty }\left(
\partial \Omega \times \left[0,T\right]\right))
\end{array}\right.
\end{equation}
where $T\in \Rset_{+}^{\ast}$ and $a_{0}\left(x,t\right)\geq 0$
in $Q$. Then for every $k\in \Nset$ there exists $P_{k}\in
\Rset\left[x\right]$ with coefficients independent of $a_{0},~f,~g$
and $h$, such that 
$$\left\Vert v\right\Vert _{\mathcal{C}^{k}\left( \overline{Q}\right) }\leq \left(\left\Vert f\right\Vert _{\mathcal{C}^{2k+1}\left( 
\overline{Q}\right) }+\left\Vert g\right\Vert _{\mathcal{C}^{2k}\left( 
\overline{\Omega }\right) }+\left\Vert h\right\Vert _{\mathcal{C}
^{2k+1}\left( \partial \Omega \times \left[ 0,T\right] \right) }\right)
P_{k}\left(\left\Vert a_{0}\right\Vert _{\mathcal{C}^{2k+1}\left( \overline{Q}
\right)}\right).$$}}
\end{Lem}

\begin{The}([\cite{cl}, Theorem 2])\label{ibvp-3.4} {\tmsamp{The solution of (\ref{ibvp-3.2.1}) is unique.}}
\end{The}

Now we are going to extend Theorems \ref{ibvp-3.2} and
\ref{ibvp-3.4} using Theorem \ref{complet-1.12}, Proposition 
\ref{complet-1.17}, Lemma \ref{complet-1.18} and the results of
section \ref{auxiliary}. The main result of this paper is the following:

\begin{The}\label{ibvp-3.5} {\tmsamp{Consider the IBVP (\ref{ibvp-3.2.1}) where
$u_{0}\in \mathcal{G}\left(\Omega\right)$. Then, there is a unique 
$u\in \mathcal{G}\left( \overline{Q}\right)$ 
such that $u$ is solution of (\ref{ibvp-3.2.1}) (with initial data $u_{0}$).}}
\end{The}
\begin{proof} Since $u_{0}\in\mathcal{G}\left(\Omega\right)$. from
Proposition \ref{complet-1.17} it follows that there exists
a sequence $\left(u_{0p}\right)_{p\geq 1}$ in
$\mathcal{G}_{c}\left(\Omega\right)$ such that 
$$u_{0p}\underset{p\to\infty}{\longrightarrow} u_{0}, ~\mbox{in}~
\left(\mathcal{G}\left(\Omega\right),\mathcal{T}_{\Omega }\right).$$
Therefore, from Theorem \ref{ibvp-3.2}, for each $p\geq 1$ there
exists $u_{p}\in\mathcal{G}\left(\overline{Q}\right)$ such
that $u_{p}$ is a solution of (\ref{ibvp-3.2.1}) with the
corresponding initial data $u_{0p}.$ Moreover, from the proof of
Theorem \ref{ibvp-3.2} it follows that $u_{0p}$ has a
representative $\widehat{u}_{0p}$ such that
$\widehat{u}_{0p}\left(\psi,-\right)\in
\mathcal{D}\left(\Omega\right),~
\forall ~\psi\in\A_{0}\left(m\right), ~p\geq 1$
and $u_{p}$ has a representative $\widehat{u}_{p}$ verifying
$\widehat{u}_{p}\left(\varphi ,-\right)\in\mathcal{C}^{\infty
}\left(\overline{Q}\right), ~\forall ~\varphi \in
\A_{0}\left(m+1\right), ~p\geq 1$
such that $\widehat{u}_{p}\left(\varphi ,-\right)$ is a classical
solution of (\ref{ibvp-3.1.1}) with initial data
$\widehat{u}_{0p}\left(\text{\textsc{I}}_{m}^{m+1}\left(\varphi\right)
,-\right)$ [indeed, from Lemma 2.1, the second
equality of (\ref{ibvp-3.1.1}): $$(*) \quad u_{p}|_{\overline{\Omega }\times
  \left\{0\right\}}=u_{0p}$$ is intended as follows: a representative of 
$u_{p}|_{\overline{\Omega }\times\left\{0\right\}}$  (resp. $u_{0p}$) is
$$R_{1}:\left(\psi,x\right)\in \A_{0}\left(m\right)\times 
\overline{\Omega }\mapsto
\widehat{u}_{p}\left(\text{\textsc{I}}_{m+1}^{m}\left(\psi\right),x,0\right)
\in  \Rset ~(\mbox{resp. $\widehat{u}_{0p}$})$$
hence, at the level of representatives, (*) can be written as
$$\widehat{u}_{p}\left(\text{\textsc{I}}_{m+1}^{m}\left(\psi\right),x,0\right)
=\widehat{u}_{0p}\left(\psi,x\right), ~\forall~\psi \in
\A_{0}\left(m\right), x\in \overline{\Omega}$$
or equivalently:
$$(**)\quad \widehat{u}_{p}\left(\varphi,x,0\right) =\widehat{u}_{0p}\left( 
\text{\textsc{I}}_{m}^{m+1}\left(\varphi\right),x\right), ~\forall ~\varphi\in \A_{0}\left(m+1\right),~\forall~x\in\overline{\Omega}.\mbox{]}$$
Therefore, from Lemma \ref{ibvp-3.1}, it follows that for each
$k\in \Nset $ there is $P_{k}\in \Rset\left[x\right]$ such that the
inequality of Lemma \ref{ibvp-3.1} holds for all $p\geq 1$:
\begin{equation}\label{ibvp-3.5.1}
\left\Vert\widehat{u}_{p}\left(\varphi,-\right)
\right\Vert_{\mathcal{C}^{k}\left(\overline{Q}\right)}\leq
P_{k}(\left\Vert\widehat{u}_{0p}\left(\text{\textsc{I}}_{m}^{m+1}\left(\varphi\right),-\right)\right\Vert_{\mathcal{C}^{2k+1}\left(\overline{\Omega}\right)},
~\forall ~\varphi \in \A_{0}\left(m+1\right).
\end{equation}
Now, we have a sequence $\left(u_{p}\right)_{p\geq 1}$ in
$\mathcal{G}\left(\overline{Q}\right)$ of solutions of
(\ref{ibvp-3.2.1}) corresponding to the initial data
$\left(u_{0p}\right)_{p\geq 1}$. We are going to prove that
the sequence $\left(u_{p}\right)_{p\geq 1}$ converges to a
$u\in\mathcal{G}\left(\overline{Q}\right)$ which is the solution
of our problem. Since $Q$ is a bounded open subset of $\Rset^{m+1}$, 
from Theorem \ref{complet-1.12} it follows that 
$\left(\mathcal{G}\left(\overline{Q}\right),\mathcal{T}_{\overline{Q},b}\right)$ is complete and hence it suffices to
show that $\left(u_{p}\right)_{p\geq 1}$ is a Cauchy sequence. 
In the sequel we  need the following notation:
$$u_{pq}:=u_{p}-u_{q} \qquad \mbox{and}\qquad
\widehat{u}_{pq}:=\widehat{u}_{p}-\widehat{u}_{q}, ~\left(p,q\in \Nset\right).$$
The definition of $u_{p}$ implies that $u_{pq}$ is a solution of
the following IBVP:
\begin{displaymath}
\left\{ 
\begin{array}{ll}
\left(u_{pq}\right)_t-\Delta u_{pq}+a_{pq}u_{pq}=0~\mbox{in}~
\mathcal{G}\left(Q\right)\\ 
u_{pq}|_{\overline{\Omega}\times\left\{0\right\}}=u_{0p}-u_{0q}
~\mbox{in}~\mathcal{G}\left(\overline{\Omega}\right)\\ 
u_{pq}|_{\partial\Omega\times\left[0,T\right]} =0 ~\mbox{in}~
\mathcal{G}\left(\partial\Omega\times\left[0,T\right]\right)
\end{array}\right. 
\end{displaymath}
where $a_{pq}:=u_{p}^{2}+u_{q}^{2}+u_{p}u_{q}.$ The definitions of 
$\widehat{u}_{p}$ and $\widehat{u}_{0p}$ shows that
$\widehat{u}_{p}\left(\varphi,-\right)$ is a classical solution to
the IBVP ~(for each $\varphi\in\A_{0}\left(m+1\right)$): 
\begin{equation}\label{ibvp-***}
\left\{ 
\begin{array}{ll}
\left(\widehat{u}_{pq}\left(\varphi,-\right)\right)_t-\Delta\widehat{
u}_{pq}\left(\varphi,-\right)+\widehat{a}_{pq}\left(\varphi,-\right) 
\widehat{u}_{pq}\left(\varphi,-\right)=0, ~\mbox{in}~\mathcal{C}^{\infty
}\left(Q\right)\\ 
\widehat{u}_{pq}\left(\varphi,-\right)|_{\overline{\Omega}\times\left\{0\right\}}=\left(\widehat{u}_{0p}-\widehat{u}_{0q}\right)\left(\text{\textsc{I}}_{m}^{m+1}\left(\varphi\right),-\right),~\mbox{in}
~\mathcal{D}\left(\Omega\right) \\ 
\widehat{u}_{pq}\left(\varphi,-\right)|_{\partial\Omega\times\left[0,T\right]}
=0,~\mbox{in} ~\mathcal{C}^{\infty }\left(\partial\Omega\times\left[0,T\right]\right)
\end{array}\right. 
\end{equation}
where
$\widehat{a}_{pq}:=\widehat{u}_{p}^{2}+\widehat{u}_{q}^{2}+\widehat{u}_{p}\widehat{u}_{q}=\left(\widehat{u}_{p}+\frac{1}{2}\widehat{u}_{q}\right)^{2}+\frac{3}{4}\widehat{u}_{q}^{2}\geq
0.$ Next, by applying Lemma \ref{ibvp-3.3} to the solution 
$\widehat{u}_{pq}\left(\varphi,-\right)$ of (\ref{ibvp-***}), with
$f=h=0$ and
$g\left(x\right):=\left(\widehat{u}_{0p}-\widehat{u}_{0q}\right)\left(\text{\textsc{I}}_{m}^{m+1}\left(\varphi\right),x\right)$,
we can conclude that for each $k\in\Nset$ there is $P_{k}\in
\Rset\left[x\right]$ (with coefficients independent  of $\widehat{a}_{pq}$ and $\widehat{u}_{0p}-\widehat{u}_{0q}$) such that (for
each $\varphi\in\A_{0}\left(m+1\right)$):
\begin{equation}\label{ibvp-3.5.2}
\left\Vert \widehat{u}_{pq}\left(\varphi,-\right)\right\Vert_{\mathcal{C}^{k}\left(\overline{Q}\right)}\leq\left\Vert\left(\widehat{u}_{0p}-\widehat{u}_{0q}\right)\left(\text{\textsc{I}}_{m}^{m+1}\left(\varphi\right),-\right)\right\Vert_{\mathcal{C}^{2k}\left(\overline{\Omega }\right)}\cdot P_{k}\left(\left\Vert 
\widehat{a}_{pq}\left(\varphi,-\right)\right\Vert _{\mathcal{C}^{2k+1}\left(\overline{Q}\right)}\right).
\end{equation}
Now, the next steps of this proof are as follows: first, we shall prove
that
{\tmsamp{
\begin{equation}\label{ibvp-3.5.3}
\left( u_{p}\right) _{p\geq 1} ~\mbox{is a
bounded sequence in} ~\left(\mathcal{G}\left(\overline{Q}\right),
\mathcal{T}_{\overline{Q},b}\right),
\end{equation}}}
and, then, by applying (\ref{ibvp-3.5.2}) and
(\ref{ibvp-3.5.3}), we shall prove that
{\tmsamp{
\begin{equation}\label{ibvp-3.5.4}
\left( u_{p}\right) _{p\geq 1} ~\mbox{is a
Cauchy sequence in}~ \left(\mathcal{G}\left(\overline{Q}\right)\mathcal{T}_{\overline{Q},b}\right),
\end{equation}}}
from which our \textit{existence} result follows. Indeed, from
(\ref{ibvp-3.5.4}) and Theorem \ref{complet-1.12} we get
$$\exists ~ u:=\underset{p\rightarrow\infty}{\lim }u_{p}\in\mathcal{G}\left(\overline{Q}\right).$$
Since from the definitions of $u_{p}$ and $u_{0p}$ we have
\begin{equation}\label{ibvp-3.5.5}
\left\{\begin{array}{ll}
\left(u_{p}\right)_{t}-\Delta u_{p}+u_{p}^{3}=0 ~\mbox{in}~\mathcal{G}
\left(Q\right)\\ 
u_{p}|_{\overline{\Omega }\times\left\{0\right\}} =u_{0p}~\mbox{in} ~\mathcal{G}\left(\overline{\Omega}\right)\\ 
u_{p}|_{\partial\Omega\times\left[0,T\right]} =0 ~\mbox{in} ~
\mathcal{G}\left(\partial\Omega\times\left[0,T\right]\right)
\end{array}\right. 
\end{equation}
and since the differential operator $$P:u\in\mathcal{G}\left(Q\right)\mapsto \left(u_{t}-\Delta u+u^{3}\right)\in \mathcal{G}\left(Q\right)$$
is continuous [this follows from [\cite{A-F-J}, Corollary 4.2] and
from the continuity
of the multiplication in $\left(\mathcal{G}\left(\Omega\right),\mathcal{T}_{\Omega}\right)]$ by taking limits in
the first equality of (\ref{ibvp-3.5.5}), we obtain
$$0=\underset{p\rightarrow\infty}{\lim
}P\left(u_{p}\right)=P\left(u\right) ~\mbox{in}~\mathcal{G}\left( Q\right).$$
Moreover, by taking limits for $p\rightarrow \infty$ in the two others
equalities of (\ref{ibvp-3.5.5}), from Lemma \ref{aux-2.1} we get
$$u|_{\overline{\Omega}\times \left\{0\right\}}=u_{0}
~\mbox{in}~\mathcal{G}\left(\overline{\Omega}\right)$$ 
and, from Lemma \ref{aux-2.2}, we have 
$$u|_{\partial\Omega\times\left[0,T\right]}=0$$
hence, $u$ is a solution of (\ref{ibvp-3.2.1}) with initial condition 
$u_{0}\in\mathcal{G}\left(\overline{\Omega}\right)$.
 It remains to prove (\ref{ibvp-3.5.3}) and (\ref{ibvp-3.5.4}).
\vspace{1cm}

\noindent \textbf{Proof of  (\ref{ibvp-3.5.3})}: In view of Lemma
\ref{complet-1.18} $(d)$ it suffices to prove the following statement:
{\tmsamp{
\begin{equation}\label{ibvp-3.5.6}
\forall ~\beta \in\Nset^{m+1}~\exists ~N\in\Nset~\mbox{such that}~
\left\Vert u_{p}\right\Vert_{\sigma }\leq \alpha_{-N}^{\bullet}
~\forall ~\sigma \leq\beta,~p\geq 1.
\end{equation}}} 
Fix any $\beta \in \Nset^{m+1}$ and choose $k\in\Nset$ such that $k\geq\left\vert\beta\right\vert$, then $k\geq \left\vert \sigma \right\vert$ for each
$\sigma\leq\beta$ and therefore, from (\ref{ibvp-3.5.1}) we have
\begin{equation}\label{ibvp-3.5.7}
\left( \widehat{u}_{p}\right)_{\sigma }\left( \varphi \right)
\leq P_{k}\left( \left\Vert \widehat{u}_{0p}\left(
      \text{\textsc{I}}_{m}^{m+1}\left(\varphi\right),-\right)\right\Vert_{\mathcal{C}^{2k+1}\left(\overline{\Omega}\right)}\right),~\forall~\varphi\in\A_{0}\left(m+1\right).
\end{equation}
Next, it is easy to check that the class in $\overline{\Rset}$ of the
moderated function defined by the second member of (\ref{ibvp-3.5.7}) 
is $$P_{k}\left(\underset{\left\vert\tau\right\vert\leq
    2k+1}{\sum}
\left\Vert u_{0p}\right\Vert_{\tau }\right);$$
in other words for each $p\geq 1$:
$$\cl\left(\varphi\mapsto P_{k}(\left\Vert\widehat{u}_{0p}\left(\text{
\textsc{I}}_{m}^{m+1}\left(\varphi\right),-\right)\right\Vert_{\mathcal{C}^{2k+1}\left(\overline{\Omega }\right)})\right)=P_{k}\left(\underset{
\left\vert\tau\right\vert\leq 2k+1}{\sum}\left\Vert u_{0p}\right\Vert_{\tau }\right).$$
Indeed, this follows from definitions and from the following trivial
remark:
"If $v\in\mathcal{E}_{M}\left(\Rset\right)$ and $P\in
\Rset\left[x\right]$ 
then $P\circ v\in \mathcal{E}_{M}\left(\Rset\right)$ and 
$\cl\left(P\circ v\right)=P\left(\cl\left(v\right)\right)$". 
Hence, from (\ref{ibvp-3.5.7}) and [\cite{A-F-J}, Lemma 2.1 (iii)] we get
\begin{equation}\label{ibvp-3.5.8}
\left\Vert u_{p}\right\Vert_{\sigma}\leq
P_{k}\left(\underset{\left\vert\tau\right\vert\leq
    2k+1}{\sum}\left\Vert u_{0p}\right\Vert_{\tau
  }\right),~\forall~\sigma\leq\beta,~\forall ~p\geq 1.
\end{equation}
Assume now that $P_{k}\left( x\right)=\overset{l}{\underset{i=0}{\sum
  }}c_{i}x^{i}$ and
$P_{k}^{\ast}\left(x\right):=\underset{i=0}{\overset{l}{\sum
  }}\left\vert c_{i}\right\vert x^{i}$.  Then clearly we have
\begin{equation}\label{ibvp-3.5.9}
x,y\in \overline{\Kset} ~\mbox{and}~ \left\vert
x\right\vert \leq \left\vert y\right\vert\Rightarrow \left\vert
P_{k}\left(x\right)\right\vert\leq P_{k}^{\ast}\left(\left\vert
x\right\vert\right)\leq P_{k}^{\ast}\left(\left\vert
y\right\vert\right).
\end{equation}
Now, we set $$t_{p}:=\underset{\left\vert\tau\right\vert\leq 2k+1}
{\sum }\left\Vert u_{0p}\right\Vert_{\tau},~\forall, ~p\geq 1,$$
since $u_{0p}\underset{p\to \infty}{\longrightarrow} u_{0}$ it is clear that
$$t_{p}\underset{p\rightarrow\infty}{\longrightarrow}
t:=\underset{\left\vert\tau\right\vert\leq 2k+1}{\sum}\left\Vert u_{0}\right\Vert_{\tau}$$
hence the set $T:=\left\{t_{p}|p\geq 1\right\}$ is bounded in 
$\left(\overline{\Rset},\mathcal{T}\right)$ which implies
\Big(see Lemma \ref{complet-1.18} $(c)$\Big) that there exists $L\in\Nset$ 
such that $T\subset V_{-L}$, that is, $$t_{p}=\left\vert t_{p}\right\vert\leq
\alpha_{-L}^{\bullet}, ~\forall~ p\geq 1.$$

\noindent Therefore from (\ref{ibvp-3.5.9}): 
$$\left\vert P_{k}\left(t_{p}\right)\right\vert\leq P_{k}^{\ast
}\left(\left\vert t_{p}\right\vert \right)\leq
P_{k}^{\ast}\left(\alpha_{-L}^{\bullet}\right), ~\forall ~p\geq 1.$$
From Lemma \ref{complet-1.18} $(c)$ there is $N\in \Nset$ such that $P_{k}^{\ast
}\left(\alpha_{-L}^{\bullet}\right)\leq\alpha_{-N}^{\bullet}$ hence
we can conclude that
$$\left\vert P_{k}\left(
    t_{p}\right)\right\vert\leq\alpha_{-N}^{\bullet}, ~\forall ~p\geq 1$$
therefore, from (\ref{ibvp-3.5.8}) we obtain
$$\left\Vert u_{p}\right\Vert_{\sigma}\leq\alpha_{-N}^{\bullet }, 
~\forall~ \sigma\leq\beta,~p\geq 1$$ which proves (\ref{ibvp-3.5.6})
and hence (\ref{ibvp-3.5.3}).
\vspace{1cm}

\noindent \textbf{Proof of (\ref{ibvp-3.5.4})}: Fix any $W_{\beta,s}$ with $\beta\in 
\Nset^{m+1}$ and $s\in \Rset$ (see Definition \ref{complet-1.2}) and choose
$k\in\Nset$ such that $k\geq\left\vert\beta\right\vert$. From the 
definition of $\left\Vert \cdot
\right\Vert_{\mathcal{C}^{2k+1}\left(\overline{Q}
\right)}$ and the previous remark on composition of a moderate 
function with a polynomial, we get
\begin{equation}\label{ibvp-3.5.10}
\cl\left(\varphi\mapsto P_{k}\left(\left\Vert
      \widehat{a}_{pq}\left(\varphi,-\right)\right\Vert_{\mathcal{C}^{2k+1}\left(\overline{Q}\right)}\right)\right)=P_{k}\left(\underset{\left\vert\sigma\right\vert\leq 2k+1}{\sum}\left\Vert a_{pq}\right\Vert_{\sigma}\right), ~\forall ~\left(p,q\right) \in \Nset^{2}.
\end{equation}
By applying (\ref{ibvp-3.5.6}) in the case
$\beta^{\ast}:=\left(2k+1,\dots,2k+1\right) \in \Nset^{m+1},$ since
$\left\vert\sigma\right\vert\leq 2k+1\Longrightarrow
\sigma_{i}\leq\beta_{i}^{\ast }=2k+1, ~\forall ~i=1,2,\dots,m$ and hence
$\sigma \leq \beta ^{\ast },$ we can conclude that there exists
$N\in\Nset$ verifying $$\left\Vert u_{p}\right\Vert_{\sigma }\leq
\alpha_{-N}^{\bullet}, ~\forall~\left\vert\sigma\right\vert\leq
2k+1,~p\geq 1.$$ Therefore, from Lemma \ref{complet-1.4}, by setting 
$C:=c\cdot\underset{\varkappa\leq\alpha}{\mathit{\max}}
\binom{\sigma}{\varkappa},$ where $c:=Card\left\{\varkappa|\varkappa\leq\sigma
\right\},$ we have 
$$\left\Vert u_{p}u_{q}\right\Vert\leq
C\cdot\alpha_{-2N}^{\bullet},~\forall ~\left\vert
  \sigma\right\vert\leq 2k+1,~\left(p,q\right) \in \Nset^{2}$$
hence $$\left\Vert a_{pq}\right\Vert_{\sigma }\leq
3C\cdot\alpha_{-2N}^{\bullet }, ~\forall~
\left\vert\sigma\right\vert\leq 2k+1,
\left(p,q\right) \in\Nset^{2}.$$ Then, if 
$N^{\prime}=Card\left\{\sigma|\left\vert\sigma\right\vert\leq
  2k+1\right\}$
we get $$\underset{\left\vert\sigma\right\vert\leq
2k+1}{\sum }\left\Vert a_{pq}\right\Vert_{\sigma }\leq
3N^{\prime}C\alpha_{-2N}^{\bullet},~\forall~\left(p,q\right) \in
\Nset^{2}.$$ Let $M\in \Nset$ such that $M>2N$ then $3N^{\prime
}C\alpha_{-2N}^{\bullet}\leq\alpha_{-M}^{\bullet}$. Hence the above
inequality implies $$\underset{\left\vert\sigma\right\vert\leq
  2k+1}{\sum}\left\Vert
  a_{pq}\right\Vert_{\sigma}\leq\alpha_{-M}^{\bullet},~\forall~\left(p,q\right)\in\Nset^{2}$$ and hence, from (\ref{ibvp-3.5.9}) it follows that
\begin{equation}\label{ibvp-3.5.11}
\left\vert P_{k}\left(\underset{\left\vert \sigma\right\vert\leq
      2k+1}{\sum}\left\Vert
      a_{pq}\right\Vert_{\sigma}\right)\right\vert\leq
P_{k}^{\ast}\left(\alpha_{-M}^{\bullet}\right),~\forall~
\left(p,q\right)\in\Nset^{2}.
\end{equation}
On the other hand for each $\left(p,q\right) \in\Nset^{2}$ we have
\begin{equation}\label{ibvp-3.5.12}
\underset{\left\vert\tau\right\vert\leq 2k}{\sum}\left\Vert u_{0p}-u_{0q}\right\Vert_{\tau}=\cl\left(\varphi\in \A_{0}\left(m+1\right)\mapsto\underset{\left\vert
\tau\right\vert\leq 2k}{\sum
}\left\Vert\partial^{\tau}\left(\widehat{u}_{0p}-\widehat{u}_{0q}\right)\left(\text{\textsc{I}}_{m}^{m+1}\left(\varphi\right),-\right)\right\Vert_{\overline{\Omega}}\right)
\end{equation}
and from the definition of $\left\Vert\cdot\right\Vert_{\mathcal{C}^{k}\left(\overline{Q}\right)}$ \ we obtain
\begin{equation}\label{ibvp-3.5.13}
\underset{\left\vert\sigma\right\vert\leq k}{\sum}\left\Vert
  u_{p}-u_{q}\right\Vert_{\sigma
}=\cl\left(\varphi\mapsto\left\Vert\left(\widehat{u}_{p}-\widehat{u}_{q}\right)\left(\varphi,-\right)\right\Vert_{\mathcal{C}^{k}\left(\overline{Q}\right)
  }\right)
\end{equation}
Now, from (\ref{ibvp-3.5.12}), (\ref{ibvp-3.5.13}),
(\ref{ibvp-3.5.10}) and [\cite{A-F-J}, Lemma 2.1 (iii)] 
we can conclude that (\ref{ibvp-3.5.2}) holds for classes:
$$\underset{\left\vert\sigma\right\vert\leq k}{\sum}\left\Vert
u_{p}-u_{q}\right\Vert_{\sigma}\leq\left(\underset{\left\vert\tau
\right\vert\leq 2k}{\sum}\left\Vert
u_{0p}-u_{0q}\right\Vert_{\tau}\right) 
\cdot P_{k}\left(\underset{\left\vert\sigma\right\vert\leq 2k+1}
{\sum}\left\Vert a_{pq}\right\Vert_{\sigma }\right),
~\forall~\left(p,q\right) 
\in\Nset^{2}$$ hence
$$\left\Vert
  u_{p}-u_{q}\right\Vert_{\sigma}\leq\left(\underset{\left\vert 
\tau\right\vert\leq 2k}{\sum}\left\Vert
u_{0p}-u_{0q}\right\Vert_{\tau}\right)\cdot P_{k}\left(
\underset{\left\vert\sigma^{\prime}\right\vert\leq 2k+1}{\sum
}\left\Vert a_{pq}\right\Vert_{\sigma^{\prime}}\right), ~\forall~\left\vert\sigma\right\vert\leq k,~\left(p,q\right)\in\Nset^{2}.$$
Obviously there exists $L\in\Nset$ such that
$$P_{k}^{\ast}\left(\alpha_{-M}^{\bullet}\right)=\left\vert
  P_{k}^{\ast}\left(\alpha_{-M}^{\bullet}\right)\right\vert\leq\alpha_{-L}^{\bullet}$$ and since $\sigma\leq\beta\Rightarrow\left\vert\sigma\right\vert\leq\left\vert\beta\right\vert\leq k$, we have
\begin{equation}\label{ibvp-3.5.14}
\left\Vert u_{p}-u_{q}\right\Vert_{\sigma
}\leq\left(\underset{\left\vert\tau\right\vert\leq 2k}{\sum
  }\left\Vert u_{0p}-u_{oq}\right\Vert_{\tau
  }\right)\cdot\alpha_{-L}^{\bullet},~ \forall~\sigma\leq
\beta,~\left(p,q\right)\in\Nset^{2}.
\end{equation}
Since $\left(u_{0p}\right)_{p\geq 1}$ is a Cauchy sequence there is
$\nu\in\Nset$ such that
\begin{equation}\label{ibvp-3.5.15}
p,q\geq\nu\Rightarrow\underset{\left\vert\tau\right\vert\leq
  2k}{\sum}\left\Vert u_{0p}-u_{0q}\right\Vert_{\tau}\leq\beta_{b}^{\bullet},
\end{equation}
where $\beta _{b}^{\bullet}$ is represented by $\widehat{\beta
}_{b}^{\bullet}:\psi\in \A_{0}\left(m\right)
\mapsto i\left(\psi\right)^{b}\in\Rset$ and $b$
 to be chosen conveniently. Clearly $\beta_{b}^{\bullet}$ can be
represented by $$^{\ast}\beta_{b}^{\bullet}:\varphi\in
\A_{0}\left(m+1\right) 
\mapsto
i\left(\text{\textsc{I}}_{m}^{m+1}\left(\varphi\right)\right)^{b}\in\Rset_{+},$$ 
we will use this representative for to prove that
\begin{equation}\label{ibvp-3.5.16}
b>\frac{m}{m+1}\left( L+s\right)
\Rightarrow\alpha_{-L}^{\bullet}\beta_{b}^{\bullet}\leq\alpha_{s}^{\bullet}.
\end{equation}
Indeed, this follows at once from [\cite{A-F-J}, Lemma 2.1~(i),~(*)] since for all
$\varphi \in \A_{0}\left( m+1\right)$ we have 
$$\left(\widehat{\alpha}_{s}^{\bullet}-\widehat{\alpha}_{-L}^{\bullet}\text{ }^{\ast}\beta_{b}^{\bullet} \right)\left(\varphi_{\varepsilon}\right)=\varepsilon^{s}\left(C_{1}-C_{2}\varepsilon^{\frac{m+1}{m}b-L-s}\right)
\geq 0$$ for $\varepsilon $ small enough ($C_{1},C_{2}$ are two 
positive constants). As a consequence, from 
(\ref{ibvp-3.5.14}), (\ref{ibvp-3.5.15}) and (\ref{ibvp-3.5.16}) we get
$$p,q\geq\nu\Rightarrow\left\Vert
  u_{p}-u_{q}\right\Vert_{\sigma}\leq\alpha_{-L}^{\bullet}\beta_{b}^{\bullet }\leq\alpha_{s}^{\bullet}, ~\forall~\sigma\leq\beta$$ or equivalently
$$p,q\geq\nu\Rightarrow (u_{p}-u_{q})\in W_{\beta,s}.$$
Finally, the uniqueness of the solution is obvious since this is precisely
[\cite{cl}, Theorem 2]. Indeed, in the proof of this result, the
initial data $u_{0}$ disappears and so, the compactness or not of
$\supp{\left(u_{0}\right)}$ is irrelevant. Therefore, this result
holds in our case.
\end{proof}

\vspace{1cm}

\noindent {\bf\underline{Acknowledgment:}} This paper is part of the second author's Ph.D. thesis written under supervision of the last one. He is grateful to CNPq-Brasil  for 
 financial support during the elaboration of his thesis (\cite{G}). He also gives thanks to his God for the gift of life and for His kindness toward him.


\bibliographystyle{mybibst}
\bibliography{biblio}

\end{document}